\documentclass{article}
\usepackage{amssymb, amsmath, amsthm}
\usepackage{MnSymbol}
\usepackage{esint}
\usepackage{wasysym}
\usepackage{graphicx}
\usepackage{color}
\usepackage{epsfig}
\usepackage{cite}

\theoremstyle{plain}
\newtheorem{thm}{Theorem}
\newtheorem{cor}[thm]{Corollary}
\newtheorem{lem}[thm]{Lemma}
\newtheorem{prop}[thm]{Proposition}

\theoremstyle{definition}
\newtheorem{defn}[thm]{Definition}
\theoremstyle{remark}

\newcommand{\norm}[1]{\left\Vert#1\right\Vert}

\newcommand{\Real}{\mathbb R}

\begin{document}
\bibliographystyle{acm}

\title{Young Measures Generated by Ideal Incompressible Fluid Flows}
\author{L\'{a}szl\'{o} Sz\'{e}kelyhidi Jr. and Emil Wiedemann}
\date{}
\maketitle
\begin{abstract}
In their seminal paper \cite{dipernamajda} R.~DiPerna and A.~Majda
introduced the notion of measure-valued solution for the incompressible
Euler equations in order to capture complex phenomena present in limits of approximate
solutions, such as persistence of oscillation and development of concentrations.
Furthermore, they gave several explicit examples exhibiting such phenomena.
In this paper we show that any measure-valued solution can be generated
by a sequence of exact weak solutions. In particular this gives rise to a very large, arguably
too large, set of weak solutions of the incompressible Euler equations.
\end{abstract}
% ----------------------------------------------------------------
\section{Introduction}
The incompressible Euler equations
\begin{equation*}
\begin{aligned}
\partial_tv+\operatorname{div}(v\otimes v)+\nabla p&=0\\
\operatorname{div}v&=0
\end{aligned}
\end{equation*}
describe the motion of an inviscid fluid with constant density in
$d$ dimensions, $d\geq2$. If we are given an initial velocity field
$v_0\in L^2(\Real^d)$ with $\operatorname{div}v_0=0$ weakly and a
positive time $0<T\leq\infty$, then the weak formulation of these equations
reads
\begin{equation}\label{weak}
\int_0^T\int_{\Real^d}(v\cdot\partial_t\phi+v\otimes
v:\nabla\phi)dxdt+\int_{\Real^d}v_0(x)\phi(x,0)dx=0;
\end{equation}
that is, we say that $v\in L^2_{loc}(\Real^d\times[0,T];\Real^d)$ is
a \emph{weak solution} if it is weakly divergence-free and satisfies
(\ref{weak}) for every $\phi\in
C_c^{\infty}(\Real^d\times[0,T);\Real^d)$ with
$\operatorname{div}\phi=0$.

Whereas classical solutions, if they exist, are unique in the class
of dissipative solutions (see \cite{lions} pp.~153-158) and moreover
conserve energy, it has been known since the seminal work of
Scheffer and Shnirelman that weak solutions are not unique and need
not conserve energy. In \cite{scheffer} V. Scheffer constructed a
weak solution in two dimensions with compact support in space and
time, thus disproving uniqueness even for zero initial data (see
also \cite{shnirel1} for a different proof). A.~Shnirelman
in \cite{shnirel2} later showed that there exist weak solutions with
decreasing energy. In \cite{euler1} and \cite{euler2} these results
were put in a unified framework based on convex integration and
Baire category techniques. In particular in \cite{euler2} the
authors show that various admissibility criteria, like energy
conservation or energy dissipation, are neither sufficient to
restore uniqueness nor can they provide for any regularity higher
than $L^{\infty}_tL^2_x$.

Several weaker concepts of solutions for Euler have arisen in the
literature, for example Brenier's generalised flows \cite{brenier1, brenier2}, Lions' dissipative solutions \cite{lions}, and
DiPerna-Majda's measure-valued solutions \cite{dipernamajda}. The
latter can be briefly described as follows: Given a sequence of
velocity fields $v_n(x,t)$, it is known from classical Young measure
theory (see e.g. \cite{young1, young2, ball, hunger, muller}) that there exists a subsequence (not
relabeled) and a parametrised probability measure $\nu_{x,t}$ on
$\Real^d$ such that for all \emph{bounded} test functions $f$,
\begin{equation*}
f(v_n(x,t))\stackrel{*}{\rightharpoonup}\int_{\Real^d}fd\nu_{x,t}
\end{equation*}
weakly* in $L^{\infty}$. One can interpret the measure $\nu_{x,t}$
as the probability distribution of the velocity field at the point
$x$ at time $t$ when the sequence $(v_n)$ exhibits faster and faster
oscillations as $n\rightarrow\infty$. Since we only have an $L^2$
bound on $(v_n)$, concentrations could occur for non-bounded $f$, in particular for the
energy density $f(v)=\frac{1}{2}|v|^2$.
DiPerna and Majda addressed this issue in \cite{dipernamajda},
providing a framework in which both oscillations and concentrations
can be described. To this end they introduced a generalised Young
measure and defined a measure-valued solution for Euler to be a
generalised Young measure that satisfies the Euler equations in an
average sense (see Section \ref{admissiblemvs} below). By considering
sequences of Leray solutions for the Navier-Stokes equations with
viscosities tending to zero, they show global existence
of measure-valued solutions for arbitrary initial data. In the context of the calculus of variations Alibert and
Bouchitt\'{e} later introduced a modified version of these
generalised Young measures \cite{alibert}, which we will work
with.

Our main result is the following (the relevant
definitions can be found in the next section):

\begin{thm}\label{mainthm} A Young measure $(\nu,\lambda,\nu^{\infty})$ on
$\Real^d$ with parameters in $\Real^d\times[0,T]$ is a
measure-valued solution of the Euler equations with bounded energy
if and only if there exists a sequence $(v_n)_{n\in\mathbb{N}}$ of
weak solutions to the Euler equations bounded in
$C\left([0,T];L^2_w(\Real^d;\Real^d)\right)$ which generate the
Young measure $(\nu,\lambda,\nu^{\infty})$ in the sense that
\begin{equation*}
f(v_n)dxdt\stackrel{*}{\rightharpoonup}\left(\int_{\Real^d}fd\nu\right) dxdt+\left(\int_{S^{d-1}}f^{\infty}d\nu^{\infty}\right)\lambda
\end{equation*}
in the sense of measures for every $f:\Real^d\rightarrow\Real$ that
possesses an $L^2$-recession function $f^{\infty}$.
\end{thm}

The proof relies on the techniques developed in \cite{euler1,euler2}, in particular on
the notion of \emph{subsolution}. As a byproduct of our analysis, we establish
a link between Euler subsolutions and measure-valued solutions in Section \ref{subsolutions}.

Theorem \ref{mainthm} shows that in a sense measure-valued solutions and weak
solutions are essentially the same for the incompressible Euler equations in dimension $d\geq 2$.
In other words, we see that in the
absence of any regularity, that is, on the level of $L^2$ or $L^{\infty}$ solutions, the notion of weak solution is too weak
to yield any information on the correlations of the velocities at
different space-time points. Indeed, measure-valued solutions merely describe the one-point
statistics $\nu_{x,t}$ of the velocity field in a weakly convergent sequence.

This is in contrast with weak and measure-valued solutions in other
contexts, such as hyperbolic conservation laws in one space
dimension, where the two defining aspects of a measure-valued
solution -- the microscopic nonlinearity and the macroscopic
conservation laws -- are strong enough to lead to compensated
compactness, see e.g.~\cite{diperna}. In such situations the
equations are usually complemented by a suitable entropy condition.
As is well known, for the incompressible Euler equations a possible
entropy condition is related to the kinetic energy $\frac{1}{2}\int |v|^2\,dx$.
Indeed, imposing the \emph{admissibility} condition, that the energy
should be bounded by the initial energy for all times, leads to the
weak-strong uniqueness for the Euler equations: any weak solution
with this initial data that satisfies the weak energy inequality is
a dissipative solution in the sense of P.-L.~Lions and therefore
coincides with the smooth solution as long as the latter exists.
For the admissibility condition in the context of hyperbolic conservation laws see \cite{dafermos,diperna}.
The weak-strong uniqueness for admissible measure-valued solutions
for the Euler equations (see Section \ref{admissiblemvs}) was
proved in \cite{weak-strong}. Here we prove:

\begin{thm}\label{energythm}
Suppose that $(\nu,\lambda,\nu^{\infty})$ is an
\emph{admissible} measure-valued solution with initial data $v_0\in
L^2(\Real^d;\Real^d)$ ($\operatorname{div}v_0=0$). Then the
generating sequence $(v_n)$ as in Theorem \ref{mainthm} may be
chosen such that in addition
\begin{equation*}
\norm{v_n(t=0)-v_0}_{L^2(\Real^d)}<\frac{1}{n}
\end{equation*}
and
\begin{equation*}
\sup_{t\in[0,T]}\frac{1}{2}\int_{\Real^d}|v_n(x,t)|^2dx\leq\frac{1}{2}\int_{\Real^d}|v_n(x,0)|^2dx.
\end{equation*}
\end{thm}

The following existence result can be easily deduced from the proof of this theorem and the existence of admissible measure-valued solutions for arbitrary $L^2$-initial data (cf. e.g. \cite{weak-strong}):  
\begin{cor}\label{c:density}
There exists an $L^2$-dense subset $\mathcal{E}$ of the set of solenoidal $L^2$-vectorfields on $\Real^d$ such that for every initial data in $\mathcal{E}$, there exist infinitely many admissible weak solutions of Euler.
\end{cor}
This is shown at the end of this work. Whether one can improve on the set $\mathcal{E}$ of such ``wild'' initial data, and obtain an existence result for
admissible weak solutions for a larger set of initial data, seems to be a very delicate issue.
In particular, such an initial data needs to be highly irregular, for otherwise we would contradict
the weak-strong uniqueness and classical local existence theorems (see e.g.~Section 2.3 in \cite{euler2}). Without the admissibility condition, existence of weak solutions has been shown in \cite{eulerexistence} for all initial data.

Finally, it should be mentioned that generalised Young measures are of
importance not only in fluid mechanics, where they emerged, but have
also been recognised a useful tool in the calculus of variations. In
particular, the question has been of some interest how Young
measures that arise from certain constrained sequences can be
characterised: The prototypic result is the theorem of Kinderlehrer
and Pedregal \cite{kipe} which states that a (classical) Young
measure is generated by a sequence of gradients if and only if it
satisfies a certain Jensen-type inequality. The result has been
generalised to so-called $\mathcal{A}$-free sequences
\cite{fonmul} and to generalised Young measures \cite{fonmulped, fonsecakruzik, k-r2}. Theorem \ref{mainthm} also gives
a characterisation of Young measures that are generated by a
constrained sequence (namely a sequence of Euler solutions), but it
differs from the previously known results in two important respects:
First, our problem does not fit into the $\mathcal{A}$-free
framework since the constant rank condition is not satisfied; and
second, our sequence not only satisfies a linear system of PDE's,
but in addition a nonlinear pointwise constraint. More concretely,
not only do we generate the Young measure with an $\mathcal{A}$-free sequence,
but with a sequence of exact solutions of the Euler
equations.\\

The rest of this paper is organised as follows: In Section
\ref{prelim} we recall the notion and key properties of generalised
Young measures and admissible measure-valued solutions. Section
\ref{mainpf} is devoted to the proof of Theorems \ref{mainthm} and
\ref{energythm}. It is split into several independent parts: First we
apply the results of \cite{euler2} to reduce the problem to finding
appropriate subsolutions in Section \ref{exact},
and we then use some more or less standard Young measure techniques in Section \ref{approx}
to reduce to discrete homogeneous oscillation Young measures.
In Section \ref{discretehom}, we present an explicit
construction of a generating sequence for discrete oscillation Young measures. Finally, in Section \ref{conclusion} we
complete the proofs of Theorem \ref{energythm} and Corollary \ref{c:density} using an argument from
\cite{euler2}.

\section{Preliminaries}\label{prelim}
\subsection{Basic Notation}

Given a locally compact separable metric space $X$, we denote by $C_c(X)$ the space of continuous functions with compact support and $C_0(X)$ the Banach space obtained from the completion of $C_c(X)$ with respect to the supremum norm. Using the Riesz representation theorem the space of finite Radon measures, denoted $\mathcal{M}(X)$, can be identified with the dual space of $C_0(X)$. We denote by $\mathcal{M}^+(X)$ and $\mathcal{M}^1(X)$ the subspaces
of positive finite measures and probability measures, respectively. 

For an open or closed subset $U\subseteq\Real^m$, $\mu\in \mathcal{M}^+(U)$ and an
open or closed subset $V\subseteq\Real^l$, we denote by
$L^{\infty}_w(U,\mu;\mathcal{M}^1(V))$ the space of
$\mu$-weakly*-measurable maps from $U$ into $\mathcal{M}^1(V)$. That
such a map $\nu$ is $\mu$-weakly*-measurable means that for each
bounded Borel function $f:V\rightarrow\Real$, the map
\begin{equation*}
x\mapsto\langle\nu_x,f\rangle:=\int_V f(z)d\nu_x(z)
\end{equation*}
is $\mu$-measurable. In case $\mu$ is the Lebesgue measure we omit
the specification of the measure.

We will denote by $L^2_x$ the space $L^2(\Real^d)$, by
$L^{\infty}_tL^2_x$ the space $L^{\infty}\left([0,T];L^2_x\right)$,
and by $CL^2_w$ the
space $C\left([0,T];L^2_w(\Real^d)\right)$ of functions that
are weakly continuous in time and $L^2$ in space; more precisely, it
is the space of maps $v:[0,T]\rightarrow L^2(\Real^d)$ such that the
map
\begin{equation*}
t\mapsto\int_{\Real^d}v(x,t)\phi(x)dx
\end{equation*}
is continuous for each test function $\phi \in L^{2}(\Real^d)$. 

We shall write $A:B$ for the scalar product of two matrices in
$\Real^{d\times d}$, that is, $A:B=\sum_{i,j}A_{ij}B_{ij}$, and
$v\otimes w$ for the tensor product of two vectors in $\Real^d$,
which is defined as a $(d\times d)$-matrix with entries $(v\otimes
w)_{ij}=v_iw_j$. Moreover we define for $v\in\Real^d$
\begin{equation*}
v\ocircle v:=v\otimes v-\frac{1}{d}|v|^2I_d,
\end{equation*}
where $I_d$ is the $d\times d$ identity matrix. Note that $v\ocircle
v$ is symmetric and has zero trace. The space of symmetric $(d\times
d)$-matrices is denoted by $\mathcal{S}^d$ and the space of
traceless symmetric $(d\times d)$-matrices by $\mathcal{S}_0^d$. If
$\phi:\Real^d\rightarrow\Real^{d\times d}$ is a differentiable
matrix-valued function, then $\operatorname{div}\phi$ is a vector
field defined by
$(\operatorname{div}\phi)_i=\sum_j\partial_{x_j}\phi_{ij}$.

If $f:X\rightarrow\Real$ and $g:Y\rightarrow\Real$ are maps from
some sets $X$, $Y$ into, say, $\Real$, then $f\otimes g$ is a map
$X\times Y\rightarrow\Real$ defined by $f\otimes g(x,y)=f(x)g(y)$,
whereas for two measures $\mu$ and $\nu$ living on two measurable
spaces $X$ and $Y$ respectively, $\mu\otimes\nu$ is a measure on
$X\times Y$ defined by $(\mu\otimes\nu)(A\times B)=\mu(A)\nu(B)$ for
measurable subsets $A\subseteq X$, $B\subseteq Y$.

Finally, $S^{d-1}\subset\Real^d$ is the $(d-1)$-dimensional unit sphere.

\subsection{Generalised Young Measures}\label{general}

In this section we recall the notion of generalised Young measure
as introduced in \cite{dipernamajda}, \cite{alibert}. For a more detailed and exhaustive discussion of
(generalised) Young measures, see e.g.~\cite{ball, hunger, k-r2, muller,weak-strong}.

Let $\Omega\subseteq\Real^m$ be an open or closed set, $p\in[1,\infty)$, and
$(w_n)_{n\in\mathbb{N}}$ a sequence of maps
$\Omega\rightarrow\Real^l$ bounded in $L^p(\Omega)$. We want to
study the limit behaviour of sequences of the form
$(f(w_n(y)))_{n\in\mathbb{N}}$ and, more 
generally, $(f(y,w_n(y)))_{n\in\mathbb{N}}$ for a certain class of test
functions $f$. Given $f\in C(\Omega\times \Real^l)$, its $L^p$-recession function $f^{\infty}$ is defined as
\begin{equation*}
f^{\infty}(y,z)=\lim_{y'\rightarrow y\atop{{z'\rightarrow
z}\atop{s\to\infty}}}\frac{f(y',sz')}{s^p}\,,
\end{equation*}
provided the limit exists. Observe that in this case $f^\infty$ is $p$-homogeneous, i.e. 
$f^{\infty}(y,\alpha z)=\alpha^p f^{\infty}(y,z)$ for all $\alpha\geq0$, $y\in\Omega, z\in\Real^l$. 
In this paper we consider test functions in the class
\begin{equation*}
\mathcal{F}_p:=\bigl\{f\in C(\Real^l):\,f^{\infty}\textrm{ exists and is continuous on }S^{l-1}\bigr\},
\end{equation*}
and more generally
\begin{equation*}
\mathcal{F}_p(\Omega):=\bigl\{f\in C(\Omega\times\Real^l):\,f^{\infty}\textrm{ exists and is continuous on }\Omega\times S^{l-1}\bigr\}.
\end{equation*}

Examples of functions in $\mathcal{F}_p(\Omega)$ are given by continuous
functions satisfying $|f(y,z)|\leq C(1+|z|^q)$ with $0\leq q<p$, in
which case $f^{\infty}=0$, or by continuous functions which are
$p$-homogeneous in $z$, in which case $f^{\infty}=f$. Of course,
functions in $\mathcal{F}_p(\Omega)$ always satisfy a bound $|f(y,z)|\leq
C(1+|z|^p)$ (where $C$ may depend on $y$, however).

A \emph{generalised Young measure on $\Real^l$ with parameters in
$\Omega$} is defined as a triple $(\nu,\lambda,\nu^{\infty})$
such that
\begin{equation*}
\nu\in L^{\infty}_w(\Omega;\mathcal{M}^1(\Real^l)),
\end{equation*}
\begin{equation*}
\lambda\in\mathcal{M}^+(\overline{\Omega}),
\end{equation*}
and
\begin{equation*}
\nu^{\infty}\in
L^{\infty}_w(\overline{\Omega},\lambda;\mathcal{M}^1(S^{l-1})).
\end{equation*}
Observe that $\nu$ is only defined Lebesgue-a.e. on $\Omega$
and $\nu^{\infty}$ is defined only $\lambda$-a.e. on
$\overline{\Omega}$. Classical Young measures are simply those where $\lambda=0$ (in which case $\nu^\infty$ is immaterial). In this case we simply write $\nu$ instead of a triple.

We are now able to state the following important result of Alibert
and Bouchitt\'{e}, which is a refinement of the construction in \cite{dipernamajda} (for proofs, see \cite{alibert},\cite{k-r2},\cite{weak-strong}):

\begin{thm}\label{fundamental}{\bf (Fundamental Theorem for Generalised Young Measures.)}\\
For $p\in[1,\infty)$ let $(w_n)_{n\in\mathbb{N}}$ be a sequence of
maps $\Omega\rightarrow\Real^l$ bounded in $L^p(\Omega)$. Then there
exists a subsequence (not relabeled) and a generalised Young measure
$(\nu,\lambda,\nu^{\infty})$ such that, for every
$f\in\mathcal{F}_p(\Omega)$,
\begin{equation*}
f(y,w_n(y))dy\stackrel{*}{\rightharpoonup}\langle\nu_y,f(y,\cdot)\rangle
dy+\langle\nu^{\infty}_y,f^{\infty}(y,\cdot)\rangle\lambda
\end{equation*}
in the sense of measures, where $\langle\nu_y,f(y,\cdot)\rangle=\int_{\Real^l}f(y,z)d\nu_y(z)$ and
$\langle\nu^{\infty}_y,f^{\infty}(y,\cdot)\rangle=\int_{S^{l-1}}f^{\infty}(y,z)d\nu_y^{\infty}(z)$.

Moreover, we then have that
$\int_{\Omega}\langle\nu_y,|\cdot|^p\rangle dy<\infty$.
\end{thm}

In the situation of the theorem, we say that the subsequence
$(w_n)$ \emph{generates} the Young measure
$(\nu,\lambda,\nu^{\infty})$ in $L^p(\Omega)$, and 
write
\begin{equation}\label{e:YMgeneration}
w_n\stackrel{\mathbf{Y_p}}{\longrightarrow}(\nu,\lambda,\nu^{\infty}).
\end{equation}
With the notation
\begin{equation*}
\llangle\nu,\lambda,\nu^{\infty};f\rrangle:=\int_{\Omega}\langle\nu,f\rangle
dy+\int_{\overline{\Omega}}\langle\nu^{\infty},f^{\infty}\rangle
d\lambda,
\end{equation*}
we can write this as
\begin{equation*}
\int_{\Omega}f(y,w_n(y))\,dy\to \llangle\nu,\lambda,\nu^{\infty};f\rrangle\quad\textrm{ for all }f\in\mathcal{F}_p(\Omega).
\end{equation*}
In the same manner we define convergence of generalised Young measures: 
we say that $(\nu^k,\lambda^k,\nu^{\infty,k})\stackrel{\mathbf{Y_p}}{\longrightarrow}(\nu,\lambda,\nu^{\infty})$ if 
\begin{equation}\label{e:YMconvergence}
\llangle\nu^k,\lambda^k,\nu^{\infty,k};f\rrangle\to\llangle\nu,\lambda,\nu^{\infty};f\rrangle
\,\textrm{ for all }f\in\mathcal{F}_p(\Omega).
\end{equation}
Indeed, \eqref{e:YMgeneration} is a special case of \eqref{e:YMconvergence}, since the function $w_n$ can be identified with the classical Young measure $\delta_{w_n}$.

The following proposition collects some well-known properties of
generalised Young measures. The proofs for the case $p=1$ can be found for instance in \cite{k-r2}, but can easily be modified for general $p\in [1,\infty)$.
\begin{prop}\label{YMprops}
\begin{enumerate}
\item[a)] There exists a countable set of functions $f_k=\phi_k\otimes h_k$, $k\in\mathbb{N}$, with $\phi_k\in C_c(\Omega)$, $h_k\in\mathcal{F}_p$ such that
$$
\llangle\nu,\lambda,\nu^{\infty};f_k\rrangle=\llangle\tilde{\nu},\tilde{\lambda},\tilde{\nu}^{\infty};f_k\rrangle\quad\forall\,k\,\Longrightarrow\, (\nu,\lambda,\nu^{\infty})=(\tilde{\nu},\tilde{\lambda},\tilde{\nu}^{\infty}).
$$
\item[b)] If $w_n\stackrel{\mathbf{Y_p}}{\longrightarrow}(\nu,\lambda,\nu^{\infty})$,
$\tilde{w}_n\stackrel{\mathbf{Y_p}}{\longrightarrow}(\tilde\nu,\tilde\lambda,\tilde\nu^{\infty})$ and $w_n-\tilde{w}_n\rightarrow0$ locally in measure,
then $\nu=\tilde\nu$.
\item[c)] If $w_n\stackrel{\mathbf{Y_p}}{\longrightarrow}(\nu,\lambda,\nu^{\infty})$ and 
$w_n-\tilde{w}_n\rightarrow0$ in $L^p_{loc}$, then 
$\tilde{w}_n\stackrel{\mathbf{Y_p}}{\longrightarrow}(\nu,\lambda,\nu^{\infty})$.
\item[d)] $w_n\rightarrow w$ strongly in $L^p_{loc}$ if and only if 
$w_n\stackrel{\mathbf{Y_p}}{\longrightarrow}\delta_w$.
\item[e)] Suppose $w_n\stackrel{\mathbf{Y_p}}{\longrightarrow}(\nu,\lambda,\nu^{\infty})$ and let $w\in L^p(\Omega)$. Then
$w_n+w\stackrel{\mathbf{Y_p}}{\longrightarrow}(\mathcal{T}_w\nu,\lambda,\nu^{\infty})$, where $\mathcal{T}_w\nu$ is the oscillation measure defined by
\begin{equation*}
\langle (\mathcal{T}_w\nu)_y,f\rangle :=\int_{\Real^l}f(z+w(y))\,d\nu_y(z)\quad\textrm{ for $f\in C_0(\Real^l)$, a.e.~$y\in\Omega$}
\end{equation*}
\end{enumerate}
\end{prop} 

In general the Young measure records the defect from strong convergence, this is signified by d). In our case the defect can come from oscillation, recorded by the
\emph{oscillation measure} $\nu$ or from concentration, recorded by the
\emph{concentration
measure} $\lambda$ and the \emph{concentration-angle measure} $\nu^{\infty}$.

In e) the Young measure $(\mathcal{T}_w\nu,\lambda,\nu^{\infty})$ is said to be the \emph{shift} of the Young measure $(\nu,\lambda,\nu^{\infty})$. This operation is useful in separating the microscopic oscillatory or concentration behaviour from the macroscopic coarse-grained state.

Part a) of the proposition implies that it suffices to test the convergence with functions $f\in\mathcal{F}_p$, i.e.~those which are \emph{independent} of $y\in\Omega$. 
A further consequence of part a) is that the convergence notion in \eqref{e:YMconvergence} is metrizable on bounded sets. This immediately leads to the following diagonal-sequence extraction principle, which we prefer to state explicitly as a proposition:
\begin{prop}\label{diagonalargument}
Suppose that for each $k\in\mathbb{N}$, 
$$
(\nu^{k,n},\lambda^{k,n},\nu^{\infty,k,n})\stackrel{\mathbf{Y_p}}{\longrightarrow}(\nu^{k},\lambda^{k},\nu^{\infty,k})\quad\textrm{ as }n\to\infty
$$
and moreover
$$
(\nu^{k},\lambda^{k},\nu^{\infty,k})\stackrel{\mathbf{Y_p}}{\longrightarrow}(\nu,\lambda,\nu^{\infty})\quad\textrm{ as }k\to\infty.
$$
Then there exists a sequence $n(k)\to\infty$ with $k\to\infty$ such that
$$
(\nu^{k,n(k)},\lambda^{k,n(k)},\nu^{\infty,k,n(k)})\stackrel{\mathbf{Y_p}}{\longrightarrow}(\nu,\lambda,\nu^{\infty})\quad\textrm{ as }k\to\infty.
$$ 
\end{prop}

\subsection{Measure-Valued Solutions of the Euler
Equations}\label{admissiblemvs}

A \emph{measure-valued solution} to the Euler equations is a
generalised Young measure on $\Real^d$ with parameters in
$\Real^d\times[0,T]$ which satisfies the Euler equations in an
average sense. This means that
\begin{equation}\label{e:average1}
\int_0^T\int_{\Real^d}\partial_t\phi\cdot\langle\nu,\xi\rangle+\nabla\phi:\langle\nu,\xi\otimes\xi\rangle
dxdt+\int_{\Real^d\times(0,T)}\nabla\phi:\langle\nu^{\infty},\theta\otimes\theta\rangle d\lambda=0
\end{equation}
for all $\phi\in C_c^\infty\left(\Real^d\times(0,T);\Real^d\right)$ with
$\operatorname{div}\phi=0$, and 
\begin{equation}\label{e:average2}
\int_{\Real^d}\nabla\psi\cdot \langle\nu_{x,t},\xi\rangle\,dx=0
\end{equation}
for all $\psi\in C_c^{\infty}(\Real^d)$ and for almost every $t$.
Here, the quantity
\begin{equation}\label{bary}
\bar{v}(x,t):=\langle\nu_{x,t},\xi\rangle
\end{equation}
is called the \emph{barycentre} of $\nu_{x,t}$ and signifies the coarse-grained, or macroscopic, flow. As usual, we have
written $\langle\nu,\xi\otimes\xi\rangle=\int\xi\otimes\xi\nu(d\xi)$
etc. 

In light of the energy bound for weak solutions of the Navier-Stokes equations,
it is natural to restrict attention to measure-valued solutions to the Euler equations which inherit this bound. 
\begin{prop}[\cite{weak-strong}]\label{technical}
Let $(v_n(x,t))$ be a sequence of functions
$\Real^d\times[0,T]\rightarrow\Real^l$ which is bounded in
$L^{\infty}\left([0,T];L^2\left(\Real^d\right)\right)$ and generates
a Young measure $(\nu,\lambda,\nu^{\infty})$ in
$L^2\left(\Real^d\times[0,T]\right)$. Then
\begin{equation}\label{e:energybound}
\operatorname{esssup}_t\left(\int_{\Real^d}\langle\nu_{x,t},|\cdot|^2\rangle
dx\right)<\infty,
\end{equation}
and the concentration measure $\lambda$ admits a disintegration of
the form 
\begin{equation}\label{e:disintegration}
d\lambda(x,t)=\lambda_t(dx)\otimes dt,
\end{equation} 
where
$t\mapsto\lambda_t$ is a bounded (w.r.t. the total variation norm)
measurable map from $[0,T]$ into
$\mathcal{M}^+\left(\Real^d\right)$.
\end{prop}

In particular, in this case Jensen's inequality implies that 
$\bar{v}(x,t)\in L^{\infty}_tL^2_x$. A well-known consequence of 
(\ref{e:average1}) is then that $\bar{v}$ can be
re-defined on a set of times of measure zero so that it belongs to
the space $CL^2_w$ (see Appendix A of \cite{euler2}), and therefore,
the initial average $\bar{v}(\cdot,0)$ is a well-defined $L^2$
function that is assumed in the sense that
$\bar{v}(\cdot,t)\rightharpoonup\bar{v}(\cdot,0)$ weakly in $L^2$ as
$t\rightarrow0$. 
Thus, we may write equation
\eqref{e:average1} as
\begin{equation}\label{avinitial}
\begin{aligned}
\int_0^T\int_{\Real^d}\partial_t\phi\cdot\langle\nu,\xi\rangle+\nabla\phi:\langle\nu,\xi\otimes\xi\rangle
dxdt+&\int_0^T\int_{\Real^d}\nabla\phi:\langle\nu^{\infty},\theta\otimes\theta\rangle\lambda_t(dx)dt\\
&=-\int_{\Real^d}\phi(x,0)\bar{v}(x,0)dx
\end{aligned}
\end{equation}
for all $\phi\in C_c^\infty\left(\Real^d\times[0,T);\Real^d\right)$ with
$\operatorname{div}\phi=0$.

Finally, the \emph{energy} of the measure-valued solution 
can be defined for almost every time $t$ as
\begin{equation}\label{YMenergy}
E(t):=\frac{1}{2}\int_{\Real^d}\langle\nu_{x,t},|\cdot|^2\rangle
dx+\frac{1}{2}\lambda_t(\Real^d).
\end{equation}

\begin{defn}{\bf (Measure-Valued Solutions.)}\label{admissible}
\begin{itemize}
\item[a)] A Young measure $(\nu,\lambda,\nu^{\infty})$ on $\Real^d$ with
parameters in $\Real^d\times[0,T]$ is called a \emph{measure-valued solution} of the Euler equations with barycentre
$\bar{v}:=\langle\nu,\xi\rangle$ if it satisfies \eqref{e:average1}-\eqref{e:average2}.
\item[b)] A Young measure $(\nu,\lambda,\nu^{\infty})$ on $\Real^d$ with
parameters in $\Real^d\times[0,T]$ is called an \emph{admissible
measure-valued solution} of the Euler equations with initial data
$v_0\in L^2(\Real^d)$ if it satisfies \eqref{e:energybound}-\eqref{e:disintegration},
 the equations \eqref{avinitial} and \eqref{e:average2} hold, and moreover
$$
E(t)\leq \frac{1}{2}\int_{\Real^d}|v_0(x)|^2dx\textrm{ for a.e. }t>0.
$$
\end{itemize}
\end{defn}
\begin{prop}\label{strongconv}
Let $(\nu,\lambda,\nu^{\infty})$ be an admissible measure-valued
solution of the Euler equations and $\bar{v}$ its barycentre as in
(\ref{bary}). Then
\begin{equation*}
\bar{v}(\cdot,t)\rightarrow\bar{v}(\cdot,0)=v_0
\end{equation*}
strongly in $L^2(\Real^d)$ as $t\rightarrow0$.
\end{prop}
\begin{proof}
We have already seen that $\bar{v}\in CL^2_w$ and therefore
\begin{equation*}
\liminf_{t\rightarrow0}\norm{\bar{v}(t)}_{L^2}\geq\norm{\bar{v}(0)}_{L^2}.
\end{equation*}
On the other hand,
\begin{equation*}
\begin{aligned}
\int|\bar{v}(t)|^2dx&=\int\left|\langle\nu_{x,t},\xi\rangle\right|^2dx\\
&\leq\int\langle\nu_{x,t},|\xi|^2\rangle
dx+\lambda_t(\Real^d)\\
&=2E(t)\leq\int|\bar{v}(0)|^2dx,
\end{aligned}
\end{equation*}
where we used the weak energy inequality in Definition
\ref{admissible}. Combining both inequalities yields
$\norm{\bar{v}(t)}_{L^2}\rightarrow\norm{\bar{v}(0)}_{L^2}$ as
$t\rightarrow0$, and since weak convergence together with
convergence of the norms implies strong convergence, we are done.
\end{proof}

\subsection{Subsolutions}\label{subsolutions}
We recall from \cite{euler1} that the Euler equations can be written in a way that 
separates them into a linear differential constraint and a nonlinear constitutive relation. 
\begin{lem}\label{tartar}
Let $v\in L^{\infty}\left([0,T];L_{loc}^2(\Real^d;\Real^d)\right)$,
$u\in
L^{\infty}\left([0,T];L_{loc}^1(\Real^d;\mathcal{S}_0^d)\right)$ and
$q$ be a distribution such that
\begin{equation}\label{linear}
\begin{aligned}
\partial_tv+\operatorname{div}u+\nabla q&=0\\
\operatorname{div}v&=0.
\end{aligned}
\end{equation}
If it also holds that
\begin{equation}\label{constraint}
u=v\ocircle v
\end{equation}
for almost every $(x,t)\in\Real^d\times[0,T]$, then $v$ and
$p:=q-\frac{1}{d}|v|^2$ are a weak solution to the Euler equations.
Conversely, if $(v,p)$ is a weak solution of Euler, then $(v,u,q)$
with $u:=v\ocircle v$ and $q:=p+\frac{1}{d}|v|^2$ solve
(\ref{linear}) and (\ref{constraint}).
\end{lem}
A pair $(v(x,t),u(x,t))$ for which there exists a pressure $q(x,t)$
such that (\ref{linear}) is satisfied is called a \emph{subsolution}
for the Euler equations. Thus, a subsolution is a solution precisely if a certain
pointwise nonlinear equation, namely \eqref{constraint}, holds.

\subsubsection{Measure-Valued Subsolutions}\label{mvssubsolutions}

The concept of subsolution easily leads to the corresponding measure-valued notion.
For the concentration-angle measure we define the set
$$
\mathbb{S}^{d-1}=\left\{(v,u)\in\Real^d\times\mathcal{S}_0^d:\,\frac{1}{d}|v|^2+|u|_{\infty}=1\right\},
$$
where $|u|_{\infty}$ denotes the operator norm of the matrix $u$. Notice that $(v,v\ocircle v)\in\mathbb{S}^{d-1}$ whenever $v\in S^{d-1}$.
Motivated by Lemma \ref{tartar}, consider a sequence 
$$
w_n=(v_n,u_n):[0,T]\times\Real^d\to \Real^d\times\mathcal{S}_0^d
$$ 
bounded in
$L^{\infty}\left([0,T];L^2(\Real^d)\times L^1(\Real^d)\right)$.
Analogously to Section \ref{general} we define the space of test-functions
\begin{equation*}
\mathcal{F}_{2,1}:=\Bigl\{f\in C(\Real^d\times\mathcal{S}_0^d):\,
f^{\infty}(v,u):=\lim_{{v'\to v\atop {u'\rightarrow
u\atop{s\rightarrow\infty}}}}\frac{f(sv',s^2u')}{s^2}\in C(\mathbb{S}^{d-1})\textrm{ exists}\Bigr\},
\end{equation*}
and similarly also the $(x,t)$-dependent version $\mathcal{F}_{2,1}(\Real^d\times[0,T])$.
We have the obvious analogue of the Fundamental Theorem for Young measures:
\begin{thm}\label{fundamental2}
Suppose $w_n=(v_n,u_n)$ is a sequence bounded in $L^{\infty}([0,T];L^2\times L^1(\Real^d))$, and $f\in\mathcal{F}_{2,1}(\Real^d\times[0,T])$. Then
there exists a subsequence (not relabeled) and a Young
measure $(\nu,\lambda,\nu^{\infty})$, with $\nu\in
L^{\infty}_w(\Real^d\times[0,T];\mathcal{M}^1(\Real^d\times\mathcal{S}_0^d))$,
$\lambda\in\mathcal{M}^+(\Real^d\times[0,T])$, $\nu^{\infty}\in
L^{\infty}_w\left(\Real^d\times[0,T],\lambda;\mathcal{M}^1(\mathbb{S}^{d-1})\right)$,
such that
\begin{equation*}
f(x,t;w_n(x,t))dxdt\stackrel{*}{\rightharpoonup}\langle\nu_{x,t},f(x,t;\cdot)\rangle
dxdt+\langle\nu^{\infty}_{x,t},f^{\infty}(x,t;\cdot)\rangle\lambda
\end{equation*}
in the sense of measures for all $f\in\mathcal{F}_{2,1}(\Real^d\times[0,T])$.
\end{thm}
In this case we write
$$
w_n\overset{\mathbf{Y_{2,1}}}{\longrightarrow}(\nu,\lambda,\nu^{\infty}).
$$

\begin{proof}
Consider the homeomorphism $\mathcal{S}^d_0\rightarrow\mathcal{S}^d_0$,
$u\mapsto|u|_\infty u$, with inverse $u\mapsto\frac{u}{\sqrt{|u|_\infty}}$. Given $f\in\mathcal{F}_{2,1}$ define 
$$
g(v,u):=f\left(\sqrt{d}v,|u|_\infty u\right).
$$
It is easy to see that $g\in \mathcal{F}_2$. Indeed, let $(v,u)$ such that $|v|^2+|u|_\infty^2=1$. 
Then 
$(\sqrt{d} v,|u|_\infty u)\in\mathbb{S}^{d-1}$ and
\begin{equation*}
\lim_{{(v',u')\rightarrow
(v,u)\atop{s\rightarrow\infty}}}\frac{g(sv',su')}{s^2}=\lim_{{(v',u')\rightarrow
(v,u)\atop{s\rightarrow\infty}}}\frac{f(\sqrt{d}sv',s^2u'|u'|_\infty)}{s^2}=f^{\infty}(\sqrt{d} v,u|u|_\infty).
\end{equation*}
Applying Theorem \ref{fundamental}
to $(\tfrac{v_n}{\sqrt{d}},\tfrac{u_n}{\sqrt{|u_n|_\infty}})$ in $L^2$ yields a generalised Young measure $(\tilde{\nu},\tilde{\lambda},\tilde{\nu}^{\infty})$ such that
\begin{equation*}
\begin{aligned}
f(v_n,u_n)dxdt=&g\left(\tfrac{v_n}{\sqrt{d}},\tfrac{u_n}{\sqrt{|u_n|_\infty}}\right)dxdt\\
\stackrel{*}{\rightharpoonup}&\langle\tilde{\nu},g\rangle
dxdt+\langle\tilde{\nu}^{\infty},g^{\infty}\rangle\tilde{\lambda}\\
=&\langle\nu,f\rangle
dxdt+\langle\nu^{\infty},f^{\infty}\rangle\lambda,
\end{aligned}
\end{equation*}
where $\tilde\lambda=\lambda$ and 
\begin{align*}
\int_{\Real^d\times\mathcal{S}_0^d} f(\xi,\zeta)d\nu(\xi,\zeta)&=\int_{\Real^d\times\mathcal{S}_0^d}
f\left(\sqrt{d}\xi,|\zeta|_\infty\zeta\right)d\tilde{\nu}(\xi,\zeta),\\
\int_{\mathbb{S}^{d-1}} f^{\infty}(\xi,\zeta)d\nu^{\infty}(\xi,\zeta)&=\int_{\{|\xi|^2+|\zeta|_\infty^2=1\}}
f^{\infty}\left(\sqrt{d}\xi,|\zeta|_\infty\zeta\right)d\tilde{\nu}^{\infty}(\xi,\zeta).
\end{align*}
\end{proof}

Let $(\nu,\lambda,\nu^{\infty})$ be a Young measure on $\Real^d\times\mathcal{S}_0^d$ with parameters in $\Real^d\times[0,T]$. Its \emph{barycentre} $\bar{w}=(\bar{v},\bar{u})$ is defined by
\begin{align}
\bar{v}(x,t)&:=\langle\nu_{x,t},\pi_1\rangle\label{barycentrev}\\
\bar{u}(x,t)&:=\langle\nu_{x,t},\pi_2\rangle
dxdt+\langle\nu^{\infty}_{x,t},\pi_2\rangle\lambda\label{barycentreu}
\end{align}
for a.e.~$x,t$, where $\pi_1$ and $\pi_2$ are the canonical
projections from $\Real^d\times\mathcal{S}_0^d$ onto $\Real^d$ and
$\mathcal{S}_0^d$, respectively. Note that $\bar{u}(x,t)$ is only a
measure. Such a Young measure is called a \emph{measure-valued subsolution} if $(\bar{v},\bar{u})$ is a subsolution in the sense of distributions, i.e. if it satisfies (\ref{linear}) for some distribution $q$.

\subsubsection{Energy and Admissibility}
\begin{defn}\label{d:energy}
For $(v,u)\in\Real^d\times\mathcal{S}_0^d$ we
define the \emph{generalised energy} by
\begin{equation*}
e(v,u):=\frac{d}{2}\lambda_{max}(v\otimes v-u),
\end{equation*}
where $\lambda_{max}$ denotes the largest eigenvalue.
\end{defn}
\begin{lem}[Lemma 3.2 in \cite{euler2}]\label{genenergy}
\begin{enumerate}
\item[a)] $e:\Real^d\times\mathcal{S}_0^d\rightarrow\Real$ is convex.
\item[b)] For every $(v,u)\in\Real^d\times\mathcal{S}_0^d$,
$\frac{1}{2}|v|^2\leq e(v,u)$, with equality if and only if
$u=v\ocircle v$.
\item[c)] For every $(v,u)\in\Real^d\times\mathcal{S}_0^d$,
$|u|_{\infty}\leq2\frac{d-1}{d}e(v,u)$, $|u|_{\infty}$ being the
operator norm of the matrix $u$.
\end{enumerate}
\end{lem}
In particular $e(v,u)\geq0$. Observe moreover that $e\in\mathcal{F}_{2,1}$ with $e^{\infty}=e$.
Then the \emph{energy} of a Young measure on $\Real^d\times\mathcal{S}_0^d$ is defined by
\begin{equation}\label{subenergy}
E(t)=\int_{\Real^d}\langle\nu_{x,t},e\rangle
dx+\int_{\Real^d}\langle\nu^{\infty}_{x,t},e\rangle\lambda_t(dx).
\end{equation}
If for a measure-valued subsolution $E(t)\leq\frac{1}{2}\int|\bar{v}(x,0)|^2dx$ for a.e. $t\geq0$, we call it an \emph{admissible} measure-valued subsolution. 

\subsubsection{Lifting}\label{lifting}
Finally, we ``lift'' measure-valued solutions to the space of measure-valued subsolutions, i.e. Young measures from $\Real^d$ to $\Real^d\times\mathcal{S}_0^d$. Let
$Q:\Real^d\rightarrow\Real^d\times\mathcal{S}_0^d$ be defined by
\begin{equation*}
Q(\xi)=(\xi,\xi\ocircle\xi).
\end{equation*}
It is easy to see that
$$
f\in\mathcal{F}_{2,1}\,\Rightarrow\,f\circ Q\in \mathcal{F}_2\quad\textrm{ with }(f\circ Q)^{\infty}=f^{\infty}\circ Q.
$$
Given now a Young measure
$(\nu,\lambda,\nu^{\infty})$ on
$\Real^d$, we define a Young measure
$(\tilde{\nu},\lambda,\tilde{\nu}^{\infty})$ on
$\Real^d\times\mathcal{S}_0^d$ by 
\begin{align*}
\langle\tilde{\nu}_{x,t},f\rangle&=\langle\nu_{x,t},f\circ Q\rangle\quad\textrm{ for $f\in C_0(\Real^d\times\mathcal{S}^d_0)$ for a.e.~$(x,t)$},\\
\langle\tilde{\nu}_{x,t}^{\infty},g\rangle&=\langle\nu^{\infty}_{x,t},g\circ Q\rangle\quad\textrm{ for $g\in C(\mathbb{S}^{d-1})$ for $\lambda$-a.e. $(x,t)$.}
\end{align*}
Then we have
\begin{prop}\label{identification}
Let $(\nu,\lambda,\nu^{\infty})$ be a measure-valued solution with bounded energy and
$(\tilde{\nu},\lambda,\tilde{\nu}^{\infty})$ be defined as
above. Suppose $(v_n,u_n)$ is bounded in $L^{\infty}_t(L^2_x\times L^1_x)$ and
$(v_n,u_n)\stackrel{\mathbf{Y_{2,1}}}{\longrightarrow}(\tilde{\nu},\lambda,\tilde{\nu}^{\infty})$. Then
\begin{enumerate}
\item[a)] the barycentre $(\bar{v},\bar{u})$ of $(\tilde{\nu},\lambda,\tilde{\nu}^{\infty})$
forms an Euler subsolution, i.e. satisfies (\ref{linear}) for some distribution $q$;
\item[b)] if $\tilde{E}(t)$ denotes the energy of the Young measure
$(\tilde{\nu},\lambda,\tilde{\nu}^{\infty})$ in the sense of
(\ref{subenergy}) and $E(t)$ the energy of $(\nu,\lambda,\nu)$ in
the sense of (\ref{YMenergy}), then $\tilde{E}(t)=E(t)$ for a.e.
$t$;
\item[c)] $v_n\stackrel{\mathbf{Y_2}}{\longrightarrow}(\nu,\lambda,\nu^{\infty})$;
\item[d)] $|u_n-v_n\ocircle v_n|\rightarrow0$
in $L^1_{loc}(\Real^d\times[0,T])$.
\end{enumerate}
\end{prop}
\begin{proof}

a) follows straightforwardly by the definition of
$(\tilde{\nu},\lambda,\tilde{\nu}^{\infty})$ and the fact
that $(\nu,\lambda,\nu^{\infty})$ is a solution to
\eqref{e:average1}-\eqref{e:average2}.

b) By definition of
$(\tilde{\nu},\lambda,\tilde{\nu}^{\infty})$, and applying
Lemma \ref{genenergy},
\begin{equation*}
\begin{aligned}
\tilde{E}(t)&=\int\langle\tilde{\nu}_{x,t},e\rangle
dx+\int\langle\tilde{\nu}^{\infty}_{x,t},e\rangle\lambda_t(dx)\\
&=\int\langle\nu,e(\xi,\xi\ocircle\xi)\rangle
dx+\int\langle\nu^{\infty},e\left(\xi,\xi\ocircle\xi\right)\rangle\lambda_t(dx)\\
&=\frac{1}{2}\int\langle\nu,|\xi|^2\rangle
dx+\frac{1}{2}\int\langle\nu^{\infty},|\xi|^2\rangle\lambda_t(dx)\\
&=\frac{1}{2}\int\langle\nu,|\xi|^2\rangle
dx+\frac{1}{2}\lambda_t(\Real^d)=E(t),
\end{aligned}
\end{equation*}
where we used that $\nu^{\infty}$ is supported on $S^{d-1}$.

c) Let $f\in\mathcal{F}_2$ and define $g:=f\circ\pi_1$.
Then $g\in\mathcal{F}_{2,1}$ with $g^{\infty}(\xi,\zeta)=f^{\infty}\left(\xi\right)$.
Therefore
\begin{equation*}
\begin{aligned}
f(v_n)dxdt=&g(v_n,u_n)dxdt\\
\stackrel{*}{\rightharpoonup}&\langle\tilde{\nu},g\rangle
dxdt+\langle\tilde{\nu}^{\infty},g^{\infty}\rangle\lambda\\
=&\langle\nu,f\rangle dxdt+\langle\nu^{\infty},f^{\infty}\rangle\lambda
\end{aligned}
\end{equation*}
by definition of $(\tilde{\nu},\lambda,\tilde{\nu}^{\infty})$ and since $g\circ Q=f$.

d) Note that the function
$f(\xi,\zeta)=\left|\zeta-\xi\ocircle\xi\right|$ belongs to
$\mathcal{F}_{2,1}$ with $f^{\infty}=f$. We can thus apply Theorem
\ref{fundamental2} with $f$ to obtain
\begin{equation*}
|u_n-v_n\ocircle
v_n|dxdt\stackrel{*}{\rightharpoonup}\langle\tilde{\nu},f\rangle
dxdt+\langle\tilde{\nu}^{\infty},f^{\infty}\rangle\lambda_tdt=0
\end{equation*}
because the set $\{(\xi,\xi\ocircle\xi):\xi\in\Real^d\}$ contains
the supports of $\tilde{\nu}$ and $\tilde{\nu}^{\infty}$,
respectively, and on this set, $f$ and $f^{\infty}$ vanish.
\end{proof}

\section{Proof of Theorems \ref{mainthm} and \ref{energythm}}\label{mainpf}
First of all observe that whenever a sequence
of weak Euler solutions bounded in $L^{\infty}_tL^2_x$ generates a generalised Young measure, then this
will be a measure-valued solution with bounded energy in the sense of Definition \ref{admissible} a). If the generating sequence consists of admissible weak solutions with initial data $v_0$, then the measure-valued solution will be admissible as in Definition \ref{admissible} b). This follows directly from the
Fundamental Theorem of Young measures (see also \cite{dipernamajda}, \cite{weak-strong}) as well as the discussion in Section \ref{admissiblemvs}.

Before we begin to prove the converse, we state a weaker version of
Theorem \ref{energythm} that we can prove along with Theorem
\ref{mainthm}. In Section \ref{conclusion} we then conclude from
this weaker statement the full assertion of Theorem \ref{energythm}.

\begin{prop}\label{almostenergythm}
Let $(\nu,\lambda,\nu^{\infty})$ be an
\emph{admissible} measure-valued solution with initial data $v_0\in
L^2(\Real^d)$. Then the
generating sequence $(v_n)$ as in Theorem \ref{mainthm} may be
chosen such that in addition
\begin{equation*}
\norm{v_n(t=0)-v_0}_{L^2(\Real^d)}<\frac{1}{n}
\end{equation*}
and
\begin{equation*}
\sup_{t\in[0,T]}\frac{1}{2}\int_{\Real^d}|v_n(x,t)|^2dx\leq\frac{1}{2}\int_{\Real^d}|v_0(x)|^2dx+\frac{1}{n}.
\end{equation*}
\end{prop}

We prove this Proposition in three steps: In Section \ref{exact} we use a result of \cite{euler2} to show that it suffices to generate measure-valued subsolutions by sequences of subsolutions. Section \ref{approx} adapts various well-known Young measure techniques to our framework to show that it suffices to construct generating sequences for discrete homogeneous oscillation measures. This is rather general and does not use any specific properties of the Euler equations. Finally, in Section \ref{discretehom} we show how to generate discrete homogeneous Young measures from subsolutions, where the plane wave analysis of the system (\ref{linear}) is exploited to give an explicit construction of the generating sequence.

\subsection{From Subsolutions to Exact Solutions}\label{exact} 
The goal of this section is to prove
\begin{prop}\label{claim1}
\begin{enumerate}
\item[a)]
We can generate $(\nu,\lambda,\nu^{\infty})$ as required in Theorem
\ref{mainthm} provided we can generate the lifted Young measure (see Subsection \ref{subsolutions}) 
$(\tilde{\nu},\tilde{\lambda},\tilde{\nu}^{\infty})$ in the sense of
Theorem \ref{fundamental2} by a sequence $(v_n,u_n)$ bounded in
$L^{\infty}_t(L^2_x\times L^1_x)$ with the properties
\begin{itemize}
\item $(v_n,u_n)$ are smooth in $\Real^d\times[0,T]$;
\item $(v_n,u_n)$ is a subsolution.
\end{itemize}
\item[b)] If $(\nu,\lambda,\nu^{\infty})$ is admissible with initial data $v_0$ and energy $E$, then we can generate it as required in Proposition \ref{almostenergythm} if the sequence $(v_n,u_n)$ additionally satisfies
\begin{itemize}
\item $\limsup_n\sup_t\int
e(v_n,u_n)dx\leq\operatorname{esssup}_tE(t)$;
\item $v_n(\cdot,0)\rightarrow v_0$ strongly in $L^2$.
\end{itemize}
\end{enumerate}
\end{prop}
\begin{proof}
Suppose now $(v_n,u_n)$ generates the Young measure
$(\tilde{\nu},\tilde{\lambda},\tilde{\nu}^{\infty})$ as in part a) of the proposition. We choose for each $n$ a function
$\bar{e}_n\in C\left(\Real^d\times(0,T);\Real\right)\cap
C\left([0,T];L^1(\Real^d;\Real)\right)$ such that
$\bar{e}_n>e_n:=e(v_n,u_n)$ on $\Real^d\times[0,T]$ and
\begin{equation}\label{energydiff}
\sup_t\int_{\Real^d}\left(\bar{e}_n-e_n\right)dx+\int_0^T\int_{\Real^d}\left(\bar{e}_n-e_n\right)dxdt<\frac{1}{n}\,.
\end{equation}

For fixed $n$, let $q_n$ be a pressure field such that
$(v_n,u_n,q_n)$ satisfies (\ref{linear}). Consider now vectorfields
$v\in C\left([0,T];L^2_w\right)$ that are smooth on
$\Real^d\times(0,T)$ and for which there exists a smooth matrix
field $u\in
C^{\infty}\left(\Real^d\times(0,T)\right);\mathcal{S}_0^d)$ such
that
\begin{equation}\label{linearn}
\begin{aligned}
\partial_tv+\operatorname{div}u+\nabla q_n&=0\\
\operatorname{div}v&=0,
\end{aligned}
\end{equation}
\begin{equation}\label{initial}
v(\cdot,0)=v_n(\cdot,0),
\end{equation}
\begin{equation}\label{final}
v(\cdot,T)=v_n(\cdot,T),
\end{equation}
and
\begin{equation}\label{targetenergy}
e\left(v(x,t),u(x,t)\right)<\bar{e}_n(x,t)
\end{equation}
for a.e. $x\in\Real^d$ and all $t\in(0,T)$. We then define the
function space $X_0^n$ by
\begin{equation*}
\begin{aligned}
X_0^n=&\left\{v\in C^{\infty}\left(\Real^d\times(0,T)\right)\cap
C\left([0,T];L^2_w\right):\right.\\
&\left.v\textit{   satisfies
(\ref{linearn}),(\ref{initial}),(\ref{final}),(\ref{targetenergy})}\right\}
\end{aligned}
\end{equation*}
and we denote the closure of $X_0^n$ under the
$C\left([0,T];L_w^2\right)$-topology by $X^n$. Note that $X_0^n$ is
non-empty, since $(v_n,u_n)\in X_0^n$. The following result is Proposition 3.3. in \cite{euler2} (the density of the set of solutions in $X^n$ is not explicitly stated in the Proposition, but is an immediate consequence of its proof):
\begin{thm}\label{dense}
The set of solutions $v\in X^n$ to the Euler equations with energy
density
\begin{equation*}
\frac{1}{2}|v(x,t)|^2=\bar{e}_n(x,t)
\end{equation*}
for every $t\in(0,T)$ and a.e. $x$ and pressure
\begin{equation*}p=q_n-\frac{1}{d}|v|^2\end{equation*} is dense in $X^n$ (w.r.t. the
$C\left([0,T];L_w^2\right)$-topology). In particular, there are infinitely many such solutions.
\end{thm}

Therefore, for $n\in\mathbb{N}$, we can find a sequence
$(v_n^k)_{k\in\mathbb{N}}\subset CL^2_w$ of weak solutions with $v_n^k\rightarrow
v_n$ as $k\rightarrow\infty$ in the $CL^2_w$-topology, i.e.
\begin{equation*}
\sup_{t\in[0,T]}\int(v_n^k-v_n)\cdot\phi\, dx\rightarrow0\quad\forall\,\phi\in L^2(\Real^d)\,.
\end{equation*}
Since $v_n\in C^{\infty}(\Real^d\times [0,T])$, we can then choose $k=k(n)$ so large that
\begin{equation*}
\sup_{t\in[0,T]}\left|\int v_n\cdot(v_n^k-v_n)dx\right|<\frac{1}{n}\,.
\end{equation*}
Next, since
$\frac{1}{2}|v_n^k|^2=\bar{e}_n$ for a.e.~$(x,t)$ by Theorem \ref{dense}, we have
\begin{equation*}
\left|\int\int|v_n^k|^2-|v_n|^2dxdt\right|\leq2\left|\int\int(\bar{e}_n-e_n)dxdt\right|+2\left|\int\int(e_n-\frac{1}{2}|v_n^2|)dxdt\right|,
\end{equation*}
where, by choice of $\bar{e}_n$, the first expression is less than
$\tfrac{1}{n}$. Concerning the last term we have
\begin{equation*}
\begin{aligned}
\left|\int\int(e_n-\frac{1}{2}|v_n^2|)dxdt\right|&=\left|\int\int\frac{d}{2}\lambda_{max}(v_n\otimes
v_n-u_n)-\frac{1}{2}|v_n|^2dxdt\right|\\
&=\left|\int\int\frac{d}{2}\lambda_{max}\left(v_n\ocircle
v_n-u_n+\frac{1}{d}|v_n|^2I_d\right)-\frac{1}{2}|v_n|^2dxdt\right|\\
&=\left|\int\int\frac{d}{2}\lambda_{max}(v_n\ocircle v_n-u_n)dxdt\right|\\
&\leq C\int\int|v_n\ocircle v_n-u_n|dxdt\rightarrow0
\end{aligned}
\end{equation*}
as $n\rightarrow\infty$, by Proposition
\ref{identification}. We also used in this calculation that for a
matrix $A$, $\lambda_{max}(A+\alpha I_d)=\lambda_{max}(A)+\alpha$.

Since
\begin{equation*}
\int\int|v_n^k-v_n|^2dxdt=\int\int|v_n^k|^2-|v_n|^2dxdt-2\int\int
v_n\cdot(v_n^k-v_n)dxdt,
\end{equation*}
we deduce that there exists a subsequence $v_n^{k(n)}$ of
Euler solutions such that $v_n-v_n^{k(n)}\rightarrow0$ in
$L^2_{loc}(\Real^d\times[0,T])$. Hence by Propositions
\ref{identification}c) and \ref{YMprops}c), this yields that the
sequence $(v_n^{k(n)})$ generates the Young measure
$(\nu,\lambda,\nu^{\infty})$ in $L^2$. This proves part a) of Proposition \ref{claim1}.

For part b), recall that in fact
$\frac{1}{2}|v_n^{k(n)}|^2=\bar{e}_n$ for a.e.~$x\in\Real^d$ for all $t\in(0,T)$. Therefore by (\ref{energydiff}) and the
assumption about the energy in part b) of Proposition \ref{claim1} we have
\begin{equation}\label{yetanotherestimate}
\limsup_n\sup_t\frac{1}{2}\int|v_n^{k(n)}|^2dx\leq
\operatorname{esssup}_tE(t)\leq\frac{1}{2}\int|v_0|^2dx.
\end{equation}
Since $v_n^{k(n)}(\cdot,0)=v_n(\cdot,0)$, we also get
\begin{equation*}
\norm{v_n^{k(n)}(\cdot,0)-v_0}_{L^2_x}=\norm{v_n(\cdot,0)-v_0}_{L^2_x}\rightarrow0
\end{equation*}
as $n\rightarrow\infty$, which, together with (\ref{yetanotherestimate}), completes the proof of the proposition.
\end{proof}

\subsection{Approximation of Generalised Young Measures}\label{approx}

In this section we reduce the problem of generating an arbitrary measure-valued solution $(\nu,\lambda,\nu^{\infty})$ to generating discrete homogeneous (i.e. independent of $x$ and $t$) oscillation measures, i.e. where $\lambda=0$ and 
\begin{equation}\label{e:discrete}
\nu_{x,t}=\sum_{i=1}^N\mu_i\delta_{(v_i,u_i)}
\end{equation}
for all $x,t$, with $\mu_i>0$, $\sum_{i=1}^N\mu_i=1$, and $(v_i,u_i)\in\Real^d\times\mathcal{S}_0^d$. This reduction is achieved by approximating
measure-valued solutions by a coupling of 
discrete oscillation measures and smooth subsolutions.
The latter represents the macroscopic flow whereas the former amounts to the microscopic oscillations encoded by the Young measure.
The general techniques for such an approximation are well known, see for instance \cite{k-r2}.

Given $k\in\mathbb{N}$ let $\mathcal{Q}^k$ be the collection of open cubes $Q_i^k\subset\Real^d\times(0,T)$ of sidelength $\tfrac{1}{k}$, whose vertices are neighbouring points on the lattice $\frac{1}{k}\mathbb{Z}^{d+1}$.  Our precise statement is the following:

\begin{thm}\label{YMapprox}
Let $(\nu,\lambda,\nu^{\infty})$ be a measure-valued subsolution with bounded energy $E(t)$. Then there exists a sequence of 
\begin{itemize}
\item discrete oscillation measures $\nu^k_{x,t}$ on $\Real^d\times\mathcal{S}_0^d$ with zero barycentres which are piecewise constant with respect to $\mathcal{Q}^k$,
\item smooth subsolutions $(\bar{v}^k,\bar{u}^k)$ bounded in $L^{\infty}_t(L^2_x\times L^1_x)$,
\end{itemize}
such that
\begin{equation}\label{approxconvergence}
\mathcal{T}_{(\bar{v}^k,\bar{u}^k)}\nu^k\overset{\mathbf{Y_{2,1}}}{\longrightarrow}(\nu,\lambda,\nu^{\infty})
\end{equation}   
and 
\begin{equation}\label{approxenergy}
\int_{\Real^d}\langle\mathcal{T}_{(\bar{v}^k,\bar{u}^k)}\nu^k,e\rangle dx\leq
\mathrm{esssup}_tE(t)+\tfrac{1}{k}\quad\forall\,t\in[0,T].
\end{equation}
Furthermore, if  $(\nu,\lambda,\nu^{\infty})$ is an admissible measure valued subsolution with initial data $v_0$, then in addition
\begin{equation}\label{e:approxinitial}
\norm{\bar{v}^k(t=0)-v_0}_{L^2(\Real^d)}<\frac{1}{k}\,.
\end{equation}

\end{thm}

\begin{proof}
Using Proposition \ref{diagonalargument} we can reduce the proof of the theorem to a series of approximation steps. 

First of all we show that a homogeneous Young measure can be approximated by discrete oscillation measures. More precisely, let $(\nu,\lambda,\nu^\infty)$ be homogeneous in the sense that $\nu$ and $\nu^{\infty}$ are independent of $x,t$ and 
$\lambda$ is a constant multiple of Lebesgue measure $\mathcal{L}^{d+1}$. The approximability is then equivalent to the following\\

\noindent{\bf Claim 1.} Let
$$
\nu\in \mathcal{M}^1(\Real^d\times\mathcal{S}^d_0),\quad\alpha\in[0,\infty),\quad \nu^{\infty}\in\mathcal{M}^1(\mathbb{S}^{d-1})
$$
such that $\langle\nu,e\rangle<\infty$, where $e$ is the generalised energy from Definition \ref{d:energy}, and assume that
\begin{equation}\label{e:zerobary}
\langle\nu,\pi_1\rangle=0,\quad\langle\nu,\pi_2\rangle+\alpha\langle\nu^{\infty},\pi_2\rangle=0.
\end{equation}
We claim that there exists a sequence $\nu^k\in\mathcal{M}^1(\Real^d\times\mathcal{S}^d_0)$ of \emph{discrete} probability measures of the form \eqref{e:discrete}
with zero barycentre, such that 
\begin{equation}\label{e:step1}
\langle\nu^k,f\rangle\to \langle\nu,f\rangle+\alpha\langle\nu^{\infty},f^\infty\rangle\quad\textrm{ for all }f\in\mathcal{F}_{2,1}\,.
\end{equation}

\noindent\textbf{Step 1. From classical to generalised measures.} 

Let us assume that $\nu,\nu^{\infty}$ are discrete probability measures, i.e.
$$
\nu=\sum_{i=1}^N\mu_i\delta_{(v_i,u_i)},\quad \nu^{\infty}=\sum_{i=1}^M\tau_i\delta_{(v_i^{\infty},u_i^{\infty})}
$$
with $(v_i,u_i)\in\Real^d\times\mathcal{S}_0^d$, $(v_i^{\infty},u_i^{\infty})\in \mathbb{S}^{d-1}$, such that 
\begin{equation*}
\sum_{i=1}^N\mu_iv_i=0,\hspace{0.4cm}\sum_{i=1}^N\mu_iu_i+\alpha\sum_{j=1}^M\tau_ju_j^{\infty}=0.
\end{equation*}
Define a sequence $(\nu^m)$ of probability measures by
\begin{equation*}
\nu^m=\left(1-\frac{1}{m}\right)\sum_{i=1}^N\mu_i\delta_{(v_{i},u_{i})}+\frac{1}{m}\sum_{j=1}^M\tau_j\delta_{\left(\sqrt{\alpha m}v_j^{\infty},\alpha mu_j^{\infty}\right)}\,,
\end{equation*}  
A direct calculation shows that the barycentre $(\bar{v}^m,\bar{u}^m)$ of $\nu^m$ satisfies $(\bar{v}^m,\bar{u}^m)\to 0$ as $m\to\infty$. Moreover, for any $f\in\mathcal{F}_{2,1}$ 
\begin{equation*}
\begin{split}
\langle\nu^m,f\rangle &=\left(1-\frac{1}{m}\right)\sum_{i=1}^N\mu_if(v_{i},u_{i})+\frac{1}{m}\sum_{j=1}^M\tau_jf\left(\sqrt{\alpha m}v_j^{\infty},\alpha mu_j^{\infty}\right)\\
&\overset{m\to\infty}\longrightarrow\sum_{i=1}^N\mu_if(v_i,u_i)+\alpha\sum_{j=1}^M\tau_jf^{\infty}(v_j^{\infty},u_j^{\infty})=\langle\nu,f\rangle+\alpha\langle\nu^{\infty},f^\infty\rangle\,.
\end{split}
\end{equation*}
Therefore also the shifted measure $\mathcal{T}_{-(\bar{v}^m,\bar{u}^m)}\nu^m$ satisfies
\begin{equation}\label{e:approx1}
\langle\mathcal{T}_{-(\bar{v}^m,\bar{u}^m)}\nu^m,f\rangle\rightarrow\langle\nu,f\rangle+\alpha\langle\nu^{\infty},f^\infty\rangle\quad\textrm{ for all }f\in\mathcal{F}_{2,1}\,.
\end{equation}

\bigskip

\noindent{\bf Step 2. From discrete to general measures with compact support}

More generally, assume that $\nu,\nu^{\infty}$ are probability measures with compact support.
By standard measure theory (see e.g.~\cite{bauer}, \S 30), we
find sequences of discrete measures $\nu^k$ with uniformly
compact support and $\nu^{k,\infty}$ such that
$$
\nu^k\stackrel{*}{\rightharpoonup}\nu\textrm{ in }\mathcal{M}(\Real^d\times\mathcal{S}_0^d),\quad \nu^{k,\infty}\stackrel{*}{\rightharpoonup}\nu^{\infty}\textrm{ in }\mathcal{M}(S).
$$
It follows easily that
\begin{equation}\label{e:approx2}
\langle\nu^k,f\rangle+\alpha\langle\nu^{k,\infty},f^{\infty}\rangle\to\langle\nu,f\rangle+\alpha\langle\nu^{\infty},f^{\infty}\rangle\quad\textrm{ for all $f\in\mathcal{F}_{2,1}$}\,.
\end{equation}

\bigskip

\noindent{\bf Step 3. From compact support to finite energy}

First of all note that the assumption $\langle \nu,e\rangle<\infty$ implies
$\langle\nu,f\rangle<\infty$ for any $f\in\mathcal{F}_{2,1}$. Indeed, any $f\in\mathcal{F}_{2,1}$ satisfies
a bound of the form $|f(\xi,\zeta)|\leq C(|\xi|^2+|\zeta|)$, therefore by Lemma \ref{genenergy}
we have $|f|\leq C'e$. Using
an idea from \cite{kruzikroubicek}, we may then approximate $\nu$ by
compactly supported measures in the following way:

For $\rho\in\mathbb{N}$, let
$r^{\rho}:\Real^d\times\mathcal{S}_0^d\rightarrow\Real$ be a smooth
function which is $1$ on $B_{\rho}$, zero on
$\left(\Real^d\times\mathcal{S}_0^d\right)\backslash B_{\rho+1}$ and
$0\leq r\leq1$ everywhere. Define also a number $s^{\rho}$ by
\begin{equation*}
s^{\rho}=\langle\nu,1-r^{\rho}\rangle,
\end{equation*}
which measures how much mass $\nu$ carries outside of $B_{\rho}(0)$.
We then define
\begin{equation*}
\nu^{\rho}:=r^{\rho}\nu+s^{\rho}\delta_{0},
\end{equation*}
which is a probability measure with support in $B_{\rho+1}$.
Heuristically, we obtain $\nu^{\rho}$ by cutting off $\nu$ outside
of $B_{\rho+1}$ and concentrating the remaining mass at zero. 
Since $r^{\rho}\to 1$ pointwise as $\rho\to\infty$, by dominated convergence we have
\begin{equation}\label{e:dct}
\langle\nu,(1-r^{\rho})f\rangle\rightarrow0\quad\textrm{ for all $f\in\mathcal{F}_{2,1}$},
\end{equation}
in particular $s^{\rho}\to 0$ and hence also
\begin{equation*}
\langle\nu-\nu^{\rho},f\rangle=\langle\nu,(1-r^{\rho})f\rangle-s^{\rho}f(0)\to 0\quad\textrm{ as }\rho\to\infty.
\end{equation*}

In order to keep the condition \eqref{e:zerobary}, let
$\bar{v}^\rho=\langle\nu^\rho,\pi_1\rangle$ and $\bar{u}^\rho=\langle\nu^\rho,\pi_2\rangle+\alpha\langle\nu^{\infty},\pi_2\rangle$ and
consider the shifted measure $\mathcal{T}_{-(\bar{v}^\rho,\bar{u}^\rho)}\nu^\rho$. 
Using \eqref{e:dct} we see that $(\bar{v}^\rho,\bar{u}^\rho)\to 0$ as $\rho\to\infty$, hence
\begin{equation}\label{e:approx3}
\langle\mathcal{T}_{-(\bar{v}^\rho,\bar{u}^\rho)}\nu^\rho,f\rangle\to\langle\nu,f\rangle\quad\textrm{ for all }f\in\mathcal{F}_{2,1}.
\end{equation}

\bigskip

Claim 1 then follows by choosing a diagonal sequence in the three approximations \eqref{e:approx1},\eqref{e:approx2},\eqref{e:approx3}.

\smallskip

Next, we show how to discretize a measure-valued subsolution so that Claim 1 can be applied to each homogeneous part separately. \\

\noindent{\bf Claim 2.} Let $(\nu,\lambda,\nu^\infty)$ be a measure-valued subsolution with bounded energy and barycentre $(\bar{v},\bar{u})$. Then there exists a sequence of smooth subsolutions $(\bar{v}^k,\bar{u}^k)$ and a sequence of generalised Young measures $(\nu^k,\lambda^k,\nu^{k,\infty})$ with zero barycentre which are piecewise constant with respect to $\mathcal{Q}^k$, such that
\begin{equation}\label{e:step2}
(\mathcal{T}_{(\bar{v}^k,\bar{u}^k)}\nu^k,\lambda^k,\nu^{k,\infty})\overset{\mathbf{Y_{2,1}}}{\longrightarrow}(\nu,\lambda,\nu^{\infty})
\end{equation}
and, if $E_k(t)$ denotes the generalised energy of $(\mathcal{T}_{(\bar{v}^k,\bar{u}^k)}\nu^k,\lambda^k,\nu^{k,\infty})$, then
\begin{equation}\label{e:step2energy}
\mathrm{esssup}_tE_k(t)\leq \mathrm{esssup}_tE(t)+\tfrac{1}{k}\,.
\end{equation}
Moreover, if $(\nu,\lambda,\nu^\infty)$ is an admissible measure-valued subsolution with initial data $v_0$, then in addition
\begin{equation}\label{e:step2initial}
\norm{\bar{v}^k(t=0)-v_0}_{L^2(\Real^d)}<\frac{1}{k}\,.
\end{equation}

\bigskip

\noindent{\bf Step 4. Regularizing.}

Let $\psi:\Real^d\rightarrow\Real$ be a standard mollification kernel, that is, smooth and
non-negative, supported on $B_1(0)$,
and $\int\psi dx=1$. Let furthermore $\chi:\Real\rightarrow\Real$ be another mollification kernel
with the same properties as $\psi$, but whose support is required to be contained in
$(-1,0)$. Define now
$\psi_{\epsilon}(x)=\frac{1}{\epsilon^d}\psi\left(\frac{x}{\epsilon}\right)$
and
$\chi_{\epsilon}(t)=\frac{1}{\epsilon}\chi\left(\frac{t}{\epsilon}\right)$,
so that the mass is still $1$ and the supports are in $B_{\epsilon}(0)$
and $(-\epsilon,0)$ respectively. Set
$\phi_{\epsilon}(x,t)=\psi_{\epsilon}(x)\chi_{\epsilon}(t)$.
We can now define for every $t\in[0,T-\epsilon]$ and
$x\in\Real^d$ another Young measure
$(\nu^{\epsilon},\lambda^{\epsilon},\nu^{\epsilon,\infty})$ by
\begin{align*}
\langle\nu^{\epsilon},f\rangle&=\langle\nu,f\rangle*\phi_{\epsilon}\quad\textrm{ for all }f\in C_0(\Real^d\times\mathcal{S}^d_0),\\
\lambda^{\epsilon}&=\lambda*\phi_{\epsilon},\\
\langle\nu^{\epsilon,\infty},g\rangle&=\frac{\left(\langle\nu^{\infty},g\rangle\lambda\right)*\phi_{\epsilon}}{\lambda^\epsilon}\quad\textrm{for all }g\in C(S).
\end{align*}
Observe that $|\langle\nu^{\infty},g\rangle|\leq \sup|g|$, so that with the above definition 
$\langle\nu^{\epsilon,\infty},g\rangle\in L^\infty(\Real^d\times(0,T-\epsilon))$, whereas $\langle\nu^{\epsilon},f\rangle,\,\lambda^{\epsilon}\in C^{\infty}(\Real^d\times[0,T-\epsilon])$.
Also, the mollified Young measure is only defined for $t\in[0,T-\epsilon]$. However, a simple rescaling of time $t\mapsto \tfrac{T}{T-\epsilon}t$ can then restore the original domain $t\in[0,T]$. Therefore, we may as well assume that $(\nu^{\epsilon},\lambda^{\epsilon},\nu^{\epsilon,\infty})$
is defined for $(x,t)\in \Real^d\times[0,T]$. Moreover, we have
\begin{equation}\label{e:approx4}
(\nu^{\epsilon},\lambda^{\epsilon},\nu^{\infty,\epsilon})\overset{\mathbf{Y_{2,1}}}{\longrightarrow}(\nu,\lambda,\nu^{\infty})\quad\textrm{ as $\epsilon\rightarrow0$,}
\end{equation}
since 
$\mu*\phi_\epsilon\stackrel{*}{\rightharpoonup}\mu$ in $\mathcal{M}$ as
$\epsilon\rightarrow0$ for any $\mu$.
Moreover, letting
$$
E_{\epsilon}(t)=\int_{\Real^d}\langle\nu^\epsilon,e\rangle\,dx+\int_{\Real^d}\langle\nu^{\infty,\epsilon},e\rangle\,\lambda^{\epsilon}_t(dx)\quad \textrm{ for }t\in[0,T],
$$
we easily see that 
$$
E_{\epsilon}(t)=E*\chi_{\epsilon}(t)\textrm{ for all }t\in[0,T].
$$
In particular for every $t\in[0,T]$ we have
\begin{equation}\label{e:mollifiedenergy}
E_{\epsilon}(t)=\int
E(t-s)\chi_{\epsilon}(s)ds\leq\mathrm{esssup}_tE(t)\int\chi=\mathrm{esssup}_tE(t).
\end{equation}

For the barycentre $(\bar{v}_{\epsilon},\bar{u}_{\epsilon})$ of this measure we have
\begin{equation}\label{e:bary}
\bar{v}_{\epsilon}=\bar{v}*\phi_{\epsilon},\quad
\bar{u}_{\epsilon}=\bar{u}*\phi_{\epsilon},
\end{equation} 
so the barycentre is
smooth and, by linearity, is a subsolution.

\smallskip

Finally, assume that $(\nu,\lambda,\nu^\infty)$ is an admissible measure-valued subsolution with initial data $v_0$.
We claim that in this case
\begin{equation}\label{e:mollifiedstrong}
\bar{v}_{\epsilon}(0)\to\bar{v}(0)\quad\textrm{ in }L^2(\Real^d),
\end{equation}
where we write $\bar{v}(t):=x\mapsto \bar{v}(x,t)$. Indeed, we have
\begin{align*}
\bar{v}_{\epsilon}(x,0)&=\int_0^{\epsilon}\int_{\Real^d}\bar{v}(x-y,s)\psi_{\epsilon}(y)\,dy\,\chi_{\epsilon}(-s)\,ds\\
&=\int_0^{\epsilon}[\bar{v}(s)*\psi_{\epsilon}]\,\chi_{\epsilon}(-s)\,ds,
\end{align*}
and moreover
\begin{align*}
\|\bar{v}(s)*\psi_{\epsilon}-\bar{v}(0)\|_{L^2}&\leq \|\bar{v}(s)*\psi_{\epsilon}-\bar{v}(0)*\psi_{\epsilon}\|_{L^2}+\|\bar{v}(0)*\psi_{\epsilon}-\bar{v}(0)\|_{L^2}\\
&\leq \|\bar{v}(s)-\bar{v}(0)\|_{L^2}+\|\bar{v}(0)*\psi_{\epsilon}-\bar{v}(0)\|_{L^2}.
\end{align*}
Therefore 
\begin{equation*}
\|\bar{v}_{\epsilon}(0)-\bar{v}(0)\|_{L^2}\leq \sup_{0<s<\epsilon}\|\bar{v}(s)-\bar{v}(0)\|_{L^2}+\|\bar{v}(0)*\psi_{\epsilon}-\bar{v}(0)\|_{L^2}.
\end{equation*}
The right hand side converges to zero as $\epsilon\to 0$ by Proposition \ref{strongconv}. This proves our claim \eqref{e:mollifiedstrong}.

\bigskip

\noindent\textbf{Step 5. Averaging.} 

Next, fix $\epsilon>0$ and consider the shifted regular Young measure
$$
(\tilde\nu,\tilde\lambda,\tilde\nu^{\infty}):=\left(\mathcal{T}_{-(\bar{v}_\epsilon,\bar{u}_\epsilon)}\nu^\epsilon,\lambda^\epsilon,\nu^{\epsilon,\infty}\right)
$$
with barycentre zero, together with the ``macroscopic'' state
$$
(\tilde v,\tilde u):=(\bar{v}_\epsilon,\bar{u}_\epsilon).
$$
We use the well-known technique of averaging, see also Lemma 4.22 in \cite{muller} and Proposition 7 in
\cite{k-r2}.
For $l\in\mathbb{N}$ let $\mathcal{Q}^l=\{Q_i^l\}_i$ be the collection of open cubes in $\Real^d\times(0,T)$ of sidelength $\frac{1}{l}$ 
with vertices on the lattice $\frac{1}{l}\mathbb{Z}^{d+1}$. 
We define $(\tilde\nu^l,\tilde{\lambda}^l,\tilde{\nu}^{l,\infty})$ by
\begin{align*}
\langle\tilde{\nu}^l,f\rangle&=\fint_{Q_i^l}\langle\tilde\nu_{x,t},f\rangle dxdt\quad\textrm{ for all }f\in C_0(\Real^d\times\mathcal{S}^d_0),\\
\tilde{\lambda}^l&=\frac{\tilde\lambda(Q_i^l)}{\mathcal{L}^{d+1}(Q_i^l)}\mathcal{L}^{d+1}\restriction_{Q_i^l}\,,\\
\langle\tilde{\nu}^{\infty,l},g\rangle&=\fint_{Q_i^l}\langle\tilde\nu^{\infty}_{x,t},g\rangle
d\tilde\lambda(x,t)\quad\textrm{ for all }g\in C(\mathbb{S}^{d-1})
\end{align*}
for all $(x,t)\in Q_i^l$ for every $i$, where $\fint_{Q_i^l}gd\mu:=\frac{1}{\mu(Q_i^l)}\int_{Q_i^l}gd\mu$ for any
measure $\mu$ and any $g\in L^1(Q_i^l;\mu)$.

Then $(\tilde{\nu}^l,\tilde{\lambda}^l,\tilde{\nu}^{\infty,l})$ is
homogeneous on each $Q_i^l$, and also has zero barycentre for a.e.~$(x,t)$.
Moreover, it follows from Proposition 8 in \cite{k-r2} that
\begin{equation}\label{e:approx5}
(\tilde\nu^l,\tilde\lambda^l,\tilde\nu^{\infty,l})\overset{\mathbf{Y_{2,1}}}{\longrightarrow}(\tilde\nu,\tilde\lambda,\tilde\nu^{\infty})\quad\textrm{ as $l\to\infty$}\,.
\end{equation}
Concerning the energy, we claim that
\begin{equation}\label{e:discreteenergy}
\sup_t\int_{\Real^d}\langle\mathcal{T}_{(\tilde{v},\tilde{u})}\tilde\nu^l,e\rangle\,dx\leq 
\sup_t\int_{\Real^d}\langle\mathcal{T}_{(\tilde{v},\tilde{u})}\tilde\nu,e\rangle\,dx+o(1)\quad\textrm{as $l\to\infty$.}
\end{equation}
To this end define
$$
(\tilde{v}^l,\tilde{u}^l)(x,t):=\fint_{Q_i^l}(\tilde{v},\tilde{u})\quad\textrm{ for $(x,t)\in Q_i^l$}
$$
and for a.e.~$t\in(0,T)$:
\begin{align*}
F_0(t)&=\int_{\Real^d}\langle\mathcal{T}_{(\tilde{v},\tilde{u})}\tilde\nu,e\rangle\,dx,\quad F_1(t)=\int_{\Real^d}\langle\mathcal{T}_{(\tilde{v}^l,\tilde{u}^l)}\tilde\nu,e\rangle\,dx,\\
F_2(t)&=\int_{\Real^d}\langle\mathcal{T}_{(\tilde{v}^l,\tilde{u}^l)}\tilde\nu^l,e\rangle\,dx,\quad F_3(t)=\int_{\Real^d}\langle\mathcal{T}_{(\tilde{v},\tilde{u})}\tilde\nu^l,e\rangle\,dx.
\end{align*}
By the definition of $\tilde\nu^l$ and since $(\tilde{v}^l,\tilde{u}^l)$ is constant on $Q_i^l$, we have
$\langle\mathcal{T}_{(\tilde{v}^l,\tilde{u}^l)}\tilde\nu^l,f\rangle=\fint_Q\langle\mathcal{T}_{(\tilde{v}^l,\tilde{u}^l)}\tilde\nu_{x,t},f\rangle dxdt$. Hence 
$$
\sup_tF_2(t)\leq \sup_tF_1(t).
$$
On the other hand $e$ satisfies the pointwise estimate
\begin{equation}\label{e:epointwise}
|e(\xi_1,\zeta_1)-e(\xi_2,\zeta_2)|\leq C(|\xi_1||\xi_1-\xi_2|+|\xi_1-\xi_2|^2+|\zeta_1-\zeta_2|),
\end{equation}
from which we obtain for a.e.~$t$
$$
|F_0(t)-F_1(t)|\leq C\left(\int_{\Real^d}\langle\tilde\nu,|\xi|^2\rangle+|\tilde{v}|^2dx\right)^{1/2}\left(\int_{\Real^d}|\tilde{v}^l-\tilde{v}|^2dx\right)^{1/2}+\int_{\Real^d}
|\tilde{v}^l-\tilde{v}|^2+|\tilde{u}^l-\tilde{u}|dx.
$$
Since $(\tilde{v},\tilde{u})\in L^{\infty}_t(L^2_x\times L^1_x)\cap C^{\infty}(\Real^d\times (0,T))$, the right hand side converges to zero as $l\to\infty$, uniformly in $t$. Similarly $\sup_t|F_2(t)-F_3(t)|$ can be made arbitrarily small. Combining these estimates yields
$$
\sup_tF_3(t)\leq \sup_tF_0(t)+o(1)\quad\textrm{ as }l\to\infty
$$
as claimed in \eqref{e:discreteenergy}.

Combining the approximations \eqref{e:approx4}-\eqref{e:bary} and \eqref{e:approx5} leads to the approximation as claimed in \eqref{e:step2}. The bound on the energy \eqref{e:step2energy} follows from \eqref{e:mollifiedenergy} and \eqref{e:discreteenergy}, and \eqref{e:step2initial} follows from \eqref{e:mollifiedstrong}. This proves
Claim 2.

\bigskip

\noindent{\bf Step 6. Conclusion of the argument.}

Given a measure-valued subsolution $(\nu,\lambda,\nu^{\infty})$, by Claim 2 there exists a sequence of generalised Young measures $(\nu^k,\lambda^k,\nu^{k,\infty})$ which are piecewise constant with respect to $\mathcal{Q}^k$, and there exist smooth subsolutions $(\bar{v}^k,\bar{u}^k)$ such that
$$
(\mathcal{T}_{(\bar{v}^k,\bar{u}^k)}\nu^k,\lambda^k,\nu^{\infty,k})\overset{\mathbf{Y_{2,1}}}{\longrightarrow}(\nu,\lambda,\nu^{\infty})\quad\textrm{ as }k\to\infty
$$
and \eqref{e:step2energy} holds.
Then, using Claim 1 we can approximate the homogeneous Young measure on each cube $Q_i^k\in\mathcal{Q}^k$ separately by a discrete oscillation measure as in \eqref{e:step1}.  In this way we 
obtain a sequence of piecewise constant discrete oscillation measures $\nu^{k,l}$ such that
$$
\nu^{k,l}\overset{\mathbf{Y_{2,1}}}{\longrightarrow}(\nu^k,\lambda^k,\nu^{k,\infty})\quad\textrm{ as }l\to\infty.
$$
In particular we have for any fixed $k\in\mathbb{N}$ and $j=0,1,\dots$
\begin{equation*}
\begin{aligned}
&\int_{\frac{j}{k}}^{\frac{(j+1)}{k}}\int_{\Real^d}\langle\mathcal{T}_{(\bar{v}^k,\bar{u}^k)}\nu^{k,l},e\rangle dxdt\rightarrow\\
&\int_{\frac{j}{k}}^{\frac{(j+1)}{k}}\int_{\Real^d}\langle\mathcal{T}_{(\bar{v}^k,\bar{u}^k)}\nu^k,e\rangle dxdt+\int_{\frac{j}{k}}^{\frac{(j+1)}{k}}\int_{\Real^d}\langle\nu^{k,\infty},e\rangle\lambda_t(dx)dt
\end{aligned}
\end{equation*}
as $l\rightarrow\infty$, and since $\nu^{k,l}$ and $(\nu^k,\lambda^k,\nu^{k,\infty})$ are $t$-independent in the time interval $(\frac{j}{k},\frac{j+1}{k})$ and $(\bar{v}^k,\bar{u}^k)$ is smooth, we conclude (analogously to \eqref{e:discreteenergy} and after passing to a subsequence) that for a.e. $t\in(0,T)$
\begin{equation*}
\int_{\Real^d}\langle\mathcal{T}_{(\bar{v}^k,\bar{u}^k)}\nu^{k,l},e\rangle dx\leq\mathrm{esssup}_t\int_{\Real^d}\langle\mathcal{T}_{(\bar{v}^k,\bar{u}^k)}\nu^k,e\rangle dx+\int_{\Real^d}\langle\nu^{k,\infty},e\rangle\lambda^k_t(dx)+\frac{1}{k}
\end{equation*}
for all $k$ and $l=l(k)$ large enough. Together with \eqref{e:step2energy} this implies \eqref{approxenergy}. 
If the measure-valued solution is admissible, we have \eqref{e:step2initial}, from which \eqref{e:approxinitial} follows. This concludes the proof.

\end{proof}

\subsection{Discrete Homogeneous Young Measures}\label{discretehom}

Let $Q=(0,1)^{d+1}$. In light of Theorem \ref{YMapprox} it remains to show the following:

\begin{prop}\label{discretegeneration}
Let
\begin{equation}\label{e:probability}
\nu=\sum_{i=1}^N\mu_i\delta_{w_i}
\end{equation}
be a probability measure on $\Real^d\times\mathcal{S}^d_0$ with zero barycentre.
Then there exists a sequence $w^k=(v^k,u^k)\in C_c^{\infty}(Q;\Real^d\times\mathcal{S}^d_0)$ of smooth subsolutions
such that 
\begin{equation}\label{e:uniformbd}
\|w^k\|_{L^{\infty}(Q)}\leq C\langle \nu,e\rangle
\end{equation}
for some fixed constant $C$, for any $f\in C(\Real^d\times\mathcal{S}^d_0)$
\begin{equation}\label{e:osc-convergence}
f(w^k)\overset{*}{\rightharpoonup}\langle \nu,f\rangle\quad\textrm{ in }L^{\infty}(\Real^d\times(0,T))
\end{equation}
and moreover if $f$ is convex, then 
\begin{equation}\label{convexestimate}
\limsup_{k\to\infty}\sup_{t\in[0,1]}\int_{[0,1]^d}f(w^k(x,t))dx\leq \langle\nu,f\rangle.
\end{equation} 
\end{prop} 

\medskip

Concerning the above proposition we remark that \eqref{e:osc-convergence} is the classical
Young measure convergence for bounded sequences. In particular it follows 
that $$\iint_{Q}f(w^k(x,t))dxdt\to \langle\nu,f\rangle.$$ The crucial point about estimate \eqref{convexestimate} is that
it is \emph{uniform} in $t$.

\bigskip

In order to construct 
generating sequences consisting of subsolutions, 
we use the localized plane-wave construction developed in \cite{euler1}.
\begin{prop}\label{p:planewaves}
Let $\bar{w}=(\bar{v},\bar{u})\in \Real^d\times\mathcal{S}_0^d$ with $\bar{v}\neq 0$. 
\begin{enumerate}
\item There exists $\eta\in S^{d}\setminus \{e_{d+1}\}$ such that
$$
w(x,t):=\bar{w}\,h((x,t)\cdot\eta)
$$
is a subsolution for any profile $h\in C^{\infty}(\Real)$.
\item There exists a second order homogeneous linear operator $\mathcal{L}_{\bar{w}}$ such that
$$
w:=\mathcal{L}_{\bar{w}}[\varphi]
$$
is a solution of \eqref{linear} for any $\varphi\in C^{\infty}(\Real^d\times\Real)$.
\item Moreover, if $\varphi(x,t)=H((x,t)\cdot\eta)$ for $H\in C^{\infty}(\Real)$, then 
$$
\mathcal{L}_{\bar{w}}[\varphi](x,t)=\bar{w}H''((x,t)\cdot\eta)
$$
\end{enumerate}
\end{prop}

\begin{proof}
Recall that we defined subsolutions as pairs $(v,u)$ such that there exists
a distribution $q$ so that \eqref{linear} holds. 

Therefore statement 1.~of the proposition
follows immediately from the plane-wave analysis of the Euler equations performed in 
Section 2 of \cite{euler1} and Remark 1 therein. Statements 2.~and 3.~follow
directly from Lemma 3.2 and Lemma 3.3 in \cite{euler1} using the identification
\begin{equation*}
(v,u,q)\mapsto U=\left( \begin{array}{cc} u+qI_d & v\\
v & 0 \end{array}\right)\,.
\end{equation*}
\end{proof}

For the proof of Proposition \ref{discretegeneration} we start with the following sharpened variant of Proposition 2.2 from \cite{euler1}:

\begin{prop}\label{p:N=2}
Let $\mu_i\geq 0$ and $w_i=(v_i,u_i)\in\Real^d\times\mathcal{S}^d_0$ for $i=1,2$ 
with $\mu_1+\mu_2=1$ and $\mu_1w_1+\mu_2w_2=0$. Assume also that
$v_1\neq v_2$. Moreover, let $f\in C(\Real^d\times\mathcal{S}^d_0)$ be convex. For any 
$\epsilon>0$ there exists 
$$
w\in C_c^{\infty}(Q;\Real^d\times\mathcal{S}^d_0)
$$
such that 
\begin{enumerate}
\item[(i)] $w=(v,u)$ is a subsolution;
\item[(ii)] $\|w\|_{L^\infty(Q)}\leq \max_i|w_i|+\epsilon$;
\item[(iii)] There exist open subsets $A_1,A_2\subset Q$ such that for $i=1,2$
$$
w(x,t)=w_i\textrm{ for }(x,t)\in A_i,\quad\textrm{ and }\quad |\mathcal{L}^{d+1}(A_i)-\mu_i|<\epsilon;
$$
\item[(iv)] For the convex function $f$ we have
\begin{equation*}%\label{e:2stepconvex}
\int_{[0,1]^d}f(w(x,t))\,dx\leq \mu_1 f(w_1)+\mu_2f(w_2)+\epsilon\quad\textrm{ for all $t\in[0,1]$.}
\end{equation*}
\end{enumerate}
Furthermore, concerning property (iii) we even have for $i=1,2$
\begin{equation}\label{e:2stepmeasure}
\left| \mathcal{L}^d\{x\in [0,1]^d:\,(x,t)\in A_i\}-\mu_i\right|<\epsilon\quad\textrm{ for all }t\in(\epsilon,1-\epsilon).
\end{equation}
\end{prop}

\begin{proof}
Fix $\delta>0$ small. 
Let $h:\Real\to\Real$ be the 1-periodic extension of
\begin{equation*}
h(s)=\begin{cases} \mu_1& \text{if
$s\in[0,\mu_2)$},\\
-\mu_2& \text{if $s\in[\mu_2,1)$}
\end{cases}
\end{equation*}
and let $h_{\delta}=h*\zeta_{\delta}$, where $\zeta\in C_c^{\infty}(-1,1)$ is a standard mollifying kernel. 
Since $h_\delta$ is 1-periodic with mean zero, there exists $H_\delta\in L^{\infty}(\Real)\cap C^{\infty}(\Real)$  such that $H_\delta''=h_\delta$. 

Let $\bar{w}=w_2-w_1$ and consider the wave direction $\eta\in S^d\setminus\{e_{d+1}\}$ as well as the 
operator $\mathcal{L}_{\bar{w}}$ obtained from Proposition \ref{p:planewaves}. For any $k\in\mathbb{N}$ set
$$
\varphi^k(x,t):=\frac{1}{k^2}H_\delta(k(x,t)\cdot\eta).
$$
Next, let $\phi\in C_c^{\infty}(Q)$ be a cutoff function, i.e. 
such that 
\begin{equation*}
0\leq \phi\leq 1\textrm{ in }Q\,\textrm{ and }\,\phi=1\textrm{ in }[\delta,1-\delta]^{d+1},
\end{equation*}
and set $w^k:=\mathcal{L}_{\bar{w}}[\phi \varphi^k]$. Then
$$
w^k\in C^{\infty}_c(Q;\Real^d\times\mathcal{S}^d_0)\quad\textrm{ and $(v^k,u^k)=w^k$ is a subsolution.}
$$
Furthermore, let
$$
W^k(x,t):=\bar{w}\,h(k(x,t)\cdot\eta)\,.
$$
Since $\mathcal{L}_{\bar{w}}$ is a homogeneous second order differential operator, 
$\mathcal{L}_{\bar{w}}[\phi\varphi^k]-\phi\mathcal{L}_{\bar{w}}[\varphi^k]$ can be written as
a sum of products of first order derivatives of $\phi$ with first order derivatives of $\varphi^k$ and of second derivatives 
of $\phi$ with $\varphi^k$. Therefore
\begin{equation*}
\|w^k-\phi\mathcal{L}_{\bar{w}}[\varphi^k]\|_{L^{\infty}(Q)}\leq \frac{C(\delta)}{k}\,.
\end{equation*}
Also, by Proposition \ref{p:planewaves} we have $\mathcal{L}_{\bar{w}}[\varphi^k](x,t)=\bar{w}\,h_\delta(k(x,t)\cdot\eta)$, 
whence we conclude using the definition of $h$ and the form of $\mathcal{L}_{\bar{w}}[\varphi^k]$ that 
\begin{equation}\label{e:supbound}
\|w^k\|_{L^{\infty}(Q)}\leq \max_i|w_i|+\frac{C(\delta)}{k}.
\end{equation}
Moreover, since $\phi=1$ on $[\delta,1-\delta]^{d+1}$, 
we actually have
\begin{equation*}
w^k(x,t)=\bar{w}\,h_{\delta}(k(x,t)\cdot\eta)\quad\textrm{ in }(\delta,1-\delta)^{d+1}.
\end{equation*}
Since $h=h_\delta$ on $(\delta,\mu_2-\delta)$ as well as on $(\mu_2+\delta,1-\delta)$, it follows that  
$w^k=w_i$ on $A_i^k$, defined by
\begin{align*}
A_1^k&=\bigl\{(x,t)\in(\delta,1-\delta)^{d+1}:\,k(x\cdot\eta)\in(\delta,\mu_2-\delta)+\mathbb{Z}\bigr\}\\
A_2^k&=\bigl\{(x,t)\in(\delta,1-\delta)^{d+1}:\,k(x\cdot\eta)\in(\mu_2+\delta,1-\delta)+\mathbb{Z}\bigr\}
\end{align*}
and furthermore
\begin{equation*}
\Bigl| \mathcal{L}^d\{x\in [0,1]^d:\,w^k(x,t)=w_i\}-\mu_i\Bigr|\leq C\delta\quad\textrm{ for }i=1,2
\end{equation*}
for any $t\in (\delta,1-\delta)$, where the constant is independent of $t$ and $\delta$.

Next, let us write $\eta=(\eta',\eta_{d+1})$ for $\eta'\in\Real^{d}$ and observe that by Proposition \ref{p:planewaves} and the assumption $v_1\neq v_2$ 
we have $\eta'\neq 0$. Fix $t\in[0,1]$ and let
$$
a=\eta_{d+1}\frac{\eta'}{|\eta'|^2},\quad x'=x+at
$$
so that $(x,t)\cdot\eta=x'\cdot\eta'$. We can then write for any $f\in C(\Real^d\times\mathcal{S}_0^d)$
\begin{equation*}
\begin{split}
\int_{[0,1]^d}f(\phi(x,t)\,W^k(x,t))\,dx&=\int_{\Real^d}f(\phi(x,t)\,W^k(x,t))\,dx\\
&=\int_{\Real^d}f\bigl(\bar{w}\,\phi(x'-at,t)h(kx'\cdot\eta')\bigr)\,dx'\\
&=\int_{\Real^d}f\bigl(\bar{w}\,\phi_t(x')h(kx'\cdot\eta')\bigr)\,dx',
\end{split}
\end{equation*}
where we have written $\phi_t(x')=\phi(x'-at,t)$. Now, standard Young measure theory implies
that for any fixed $t$
\begin{equation*}
\int_{\Real^d}f\bigl(\bar{w}\,\phi_t(x')h(kx'\cdot\eta')\bigr)\,dx'\overset{k\to\infty}{\longrightarrow}\int_{\Real^d}\bigl[\mu_1f(w_1\,\phi_t(x'))+\mu_2f(w_2\,\phi_t(x'))\bigr]\,dx'\,.
\end{equation*}
Moreover, since the family $\{\phi_t\}_{t\in[0,1]}$ is equicontinuous with $\|\phi_t-\phi_{t'}\|_{L^{\infty}(\Real^d)}\leq C|t-t'|$ for some fixed constant $C$, 
the convergence above in fact holds uniformly in $t\in[0,1]$. Furthermore, if $f$ is convex, then, since $|\phi_t(x')|\leq 1$ for all $x'$,
\begin{equation*}
\begin{split}
\int_{\Real^d}\bigl[\mu_1f(w_1\,\phi_t(x'))+\mu_2f(w_2\,\phi_t(x'))&\bigr]\,dx'\leq \int_{\Real^d}\phi_t(x')\bigl[\mu_1f(w_1)+\mu_2f(w_2)\bigr]\,dx'\\
&=\int_{[0,1]^d}\phi(x,t)\bigl[\mu_1f(w_1)+\mu_2f(w_2)\bigr]\,dx\\
&\leq \mu_1f(w_1)+\mu_2f(w_2).
\end{split}
\end{equation*}
 Since any continuous convex function is locally Lipschitz, the $L^\infty$ bound \eqref{e:supbound} together with the uniform $L^{\infty}$-boundedness of $\phi W^k$ and the fact that $w^k=\phi W^k$ on $A_1^k\cup A_2^k$
 implies that, by choosing first $\delta>0$ sufficiently small and then $k$ sufficiently large, we can ensure that
 $$
 \left|\int_{[0,1]^d} f(w^k(x,t))\,dx-\int_{[0,1]^d} f(\phi(x,t)\,W^k(x,t))\,dx\right|\leq \epsilon/2\textrm{ for all }t\in[0,1].
 $$
Consequently, for $k$ sufficiently large we obtain
$$
\int_{[0,1]^d}f(w^k(x,t))\,dx\leq \mu_1f(w_1)+\mu_2f(w_2)+\epsilon.
$$
This concludes the proof.

\end{proof}

\begin{proof}[Proof of Proposition \ref{discretegeneration}]
Let  
$$
\sum_{i=1}^N\mu_i\delta_{w_i}
$$
with $w_i=(v_i,u_i)\in \Real^d\times\mathcal{S}^d_0$ be a probability measure with barycentre zero, and such that
\begin{equation}\label{nondeg}
\textrm{span}(v_1,\dots,v_N)\,\textrm{ has maximal rank}
\end{equation}
subject to the constraint $\sum_i\mu_iv_i=0$, i.e. dimension$=\min(d,N-1)$.

We prove by induction on $N$ that for any $\epsilon>0$ and a given convex function $f\in C(\Real^d\times\mathcal{S}^d_0)$ there exists
a smooth subsolution $w=(v,u)\in C_c^{\infty}(Q;\Real^d\times\mathcal{S}^d_0)$ such that
\begin{enumerate}
\item[(i)] $w=(v,u)$ is a subsolution;
\item[(ii)] $\|w\|_{L^\infty(Q)}\leq \max_i|w_i|+\epsilon$;
\item[(iii)] There exist open subsets $A_i\subset Q$ for $i=1,\dots,N$ such that
$$
w(x,t)=w_i\textrm{ for }(x,t)\in A_i,\quad\textrm{ and }\quad |\mathcal{L}^{d+1}(A_i)-\mu_i|<\epsilon;
$$
\item[(iv)]  
\begin{equation*}
\int_{[0,1]^d}f(w(x,t))\,dx\leq \sum_{i=1}^N\mu_if(w_i)+\epsilon\quad\textrm{ for all $t\in[0,1]$.}
\end{equation*}
\end{enumerate}
Proposition \ref{discretegeneration} then follows easily. Indeed, the non-degeneracy assumption \eqref{nondeg} can be achieved by perturbing $v_i$ slightly,
and then (ii),(iii),(iv) implies \eqref{e:uniformbd}, \eqref{e:osc-convergence} and \eqref{convexestimate}, respectively. Observe in particular that we can find a fixed sequence $w^k$ as required in Proposition \ref{discretegeneration} so that \eqref{convexestimate} holds for \emph{all} convex $f$; indeed, this can be obtained by a diagonal argument and the observation that it suffices to show \eqref{convexestimate} for a countable set of convex functions.  

\smallskip

The case $N=1$ is trivial, since then $\nu=\delta_0$ and we can simply take $w\equiv 0$. The case $N=2$ is directly handled
by Proposition \ref{p:N=2}, since the non-degeneracy condition \eqref{nondeg} implies $v_1\neq v_2$. 

\smallskip

\noindent{\bf Induction step. } Let $\nu=\sum_{i=1}^{N+1}\mu_i\delta_{w_i}$ for $N\geq 2$. Using condition \eqref{nondeg} and by a reordering of the vectors $v_1\dots v_{N+1}$ if necessary, we may assume without loss of generality that 
$$
v_{N+1}\neq 0\textrm{ and }(v_1,\dots,v_N)\textrm{ satisfies \eqref{nondeg}}.
$$
Define 
the probability measures
\begin{align*}
\nu_1&:=\mu_{N+1}\delta_{w_{N+1}}+(1-\mu_{N+1})\delta_{\bar{w}}\\
\nu_2&:=\sum_{i=1}^N\frac{\mu_i}{1-\mu_{N+1}}\delta_{w_i-\bar{w}}, 
\end{align*}
where
\begin{equation*}
\bar{w}:=\sum_{i=1}^N\frac{\mu_i}{1-\mu_{N+1}}w_i\,.
\end{equation*}
Observe that both $\nu_1$ and $\nu_2$ have zero barycentre. Moreover, by a direct calculation we check that $v_{N+1}\neq\bar{v}$. Therefore $\nu_1$ satisfies the assumptions of Proposition \ref{p:N=2} and $\nu_2$ satisfies the induction hypothesis.
Therefore, given $\epsilon>0$ we obtain two subsolutions
$$
W_1,W_2\in C_c^\infty(Q;\Real^d\times\mathcal{S}^d_0)
$$
satisfying properties (i)-(iv) with respect to the measures $\nu_1,\nu_2$, and moreover
$W_1$ satisfies the time-slice estimate \eqref{e:2stepmeasure}.

Let $A,B\subset Q$ be the open sets from property (iii) for $W_1$, with $B$ corresponding to the value $\bar{w}$, i.e.~such that
\begin{equation}\label{e:B}
W_1(x,t)=\bar{w}\textrm{ for }(x,t)\in B,\quad\textrm{ and }\quad |\mathcal{L}^{d+1}(B)-(1-\mu_{N+1})|<\epsilon,
\end{equation}
and let 
$$
B_{\epsilon}=\{(x,t)\in B:\,\epsilon<t<1-\epsilon\}.
$$
Fix a finite family of disjoint cubes 
\begin{equation*}
\tilde B=\bigcup_{j=1}^{M}((x_j,t_j)+\alpha_j Q),\hspace{0.4cm}(x_j,t_j) \in Q,\hspace{0.3cm}\alpha_j>0,
\end{equation*}
such that 
\begin{equation}\label{cubeexhaustion}
\tilde B\subset B_{\epsilon},\qquad \mathcal{L}^{d+1}(B_\epsilon \backslash \tilde B)<\epsilon,
\end{equation}
and set
\begin{equation*}
w(x,t)=W_1(x,t)+\sum_{j=1}^{M}W_2\bigl(\tfrac{x-x_j}{\alpha_j},\tfrac{t-t_j}{\alpha_j}\bigr).
\end{equation*}
We claim that $w$ satisfies (i)-(iv) for the measure $\nu$.

To start with, it is easy to see by linearity that $w\in C_c^\infty(Q;\Real^d\times\mathcal{S}^d_0)$ and (i),(ii) hold.
Concerning (iii), let  $A_{N+1}=A$ and $A_i\subset Q$ for $i=1,\dots,N$ be defined by 
$$
A_i=\bigcup_{j=1}^{M}((x_j,t_j)+\alpha_j \tilde A_i),
$$
where $\tilde A_i$ are the open sets from property (iii) for $W_2$.
Then
$$
w(x,t)=w_i\textrm{ for }(x,t)\in A_i\,\textrm{ for }i=1,\dots,N+1,
$$
and moreover, for $i=1,\dots,N$
$$
\mathcal{L}^{d+1}(A_i)=\sum_{j=1}^M\alpha_j^{d+1}\mathcal{L}^{d+1}(\tilde A_i).
$$
On the other hand $\sum_j\alpha_j^{d+1}=\mathcal{L}^{d+1}(\tilde B)$, hence
from \eqref{cubeexhaustion} and the definition of $B_{\epsilon}$ we deduce $|\mathcal{L}^{d+1}(B)-\sum_j\alpha_j^{d+1}|<5\epsilon$. Combined with \eqref{e:B} this yields 
$\bigl|1-\mu_{N+1}-\sum_j\alpha_j^{d+1}\bigr|<6\epsilon$, therefore we obtain
\begin{equation*}
\begin{split}
\bigl|\mathcal{L}^{d+1}(A_i)-\mu_i\bigr|&\leq \bigl|\sum_j\alpha_j^{d+1}(\mathcal{L}^{d+1}(\tilde A_i)-\tfrac{\mu_i}{1-\mu_{N+1}})\bigr|+\frac{\mu_i}{1-\mu_{N+1}}\bigl|1-\mu_{N+1}-\sum_j\alpha_j^{d+1}\bigr|\\
&\leq C\epsilon
\end{split}
\end{equation*}
for some fixed constant $C$. Thus, by replacing $\epsilon$ by $\epsilon/C$ in the above we can ensure property (iii) for $w$.

Next, recall the convex function $f\in C(\Real^d\times\mathcal{S}^d_0)$. Concerning property (iv) for $w$, assume first that $t\in [0,\epsilon]\cup[1-\epsilon,1]$.
Note that in this case $w(x,t)=W_1(x,t)$ for all $x$, hence using property (iv) for $W_1$ and the convexity of $f$ we obtain
\begin{equation*}
\begin{split}
\int_{[0,1]^d}f(w(x,t))\,dx&=\int_{[0,1]^{d}}f(W_1(x,t))dx\\
&\leq \mu_{N+1} f(w_{N+1})+(1-\mu_{N+1})f(\bar{w})+\epsilon\\
&\leq \mu_{N+1} f(w_{N+1})+\sum_{i=1}^N\mu_if(w_i)+\epsilon.
\end{split}
\end{equation*} 

Next, let $t\in(\epsilon,1-\epsilon)$, and recall that by Proposition \eqref{p:N=2} the function $W_1$ satisfies additionally the time-slice estimate \eqref{e:2stepmeasure}. 
Define the sets 
\begin{equation*}
A_t=\{x:\,(x,t)\in A\},\quad B_t=\{x:\,(x,t)\in B\},\quad \tilde B_t=\{x:\,(x,t)\in \tilde B\}.
\end{equation*}
By estimate (\ref{e:2stepmeasure}) we have
\begin{equation*}
\left|\mathcal{L}^d(A_t)-\mu_{N+1}\right|<\epsilon,\quad
\left|\mathcal{L}^d(B_t)-(1-\mu_{N+1})\right|<\epsilon,
\end{equation*}
hence $\mathcal{L}^d([0,1]^d\setminus (A_t\cup B_t))<2\epsilon$.
We can therefore estimate
\begin{equation*}
\begin{aligned}
&\int_{[0,1]^d}f(w(x,t))dx\\
\leq &\int_{A_t}f(w(x,t))dx+\int_{B_t\setminus \tilde B_t}f(w(x,t))dx+\int_{\tilde B_t}f(w(x,t))dx+O(\epsilon)\\
=&\mathcal{L}^d(A_t)f(w_{N+1})+\mathcal{L}^d(B_t\setminus\tilde B_t)f(\bar{w})+\int_{\tilde B_t}f(w(x,t))dx+O(\epsilon).
\end{aligned}
\end{equation*}
For the third integral we calculate
\begin{equation*}
\begin{split}
\int_{\tilde B_t}f(w(x,t))dx&=\sideset{}{'}\sum_j\int_{x_j+\alpha_j[0,1]^d}f\Bigl(\bar{w}+W_2\bigl(\tfrac{x-x_j}{\alpha_j},\tfrac{t-t_j}{\alpha_j}\bigr)\Bigr)\,dx\\
&=\sideset{}{'}\sum_j\alpha_j^d\int_{[0,1]^d}f\Bigl(\bar{w}+W_2\bigl(x,\tfrac{t-t_j}{\alpha_j}\bigr)\Bigr)\,dx\\
&\leq \mathcal{L}^d(\tilde B_t)\sum_{i=1}^N\frac{\mu_i}{1-\mu_{N+1}}f(w_i)+\epsilon,
\end{split}
\end{equation*}
where $\sideset{}{'}\sum_j$ denotes summation over those $j$ for which $t \in t_j+(0,\alpha_j)$, and in the last inequality we have used property (iv) for $W_2$. 
Furthermore, by convexity of $f$ we also have $f(\bar{w})\leq \sum_{i=1}^N\frac{\mu_i}{1-\mu_{N+1}}f(w_i)$. 
It follows that
\begin{equation*}
\int_{[0,1]^d}f(w(x,t))dx\leq \sum_{i=1}^{N+1}\mu_if(w_i)+O(\epsilon),
\end{equation*}
implying property (iv) for $w$ as required.

\end{proof}

\subsection{Proof of Proposition \ref{almostenergythm} and Theorem \ref{energythm}}\label{conclusion} 

\subsubsection{Proof of Proposition \ref{almostenergythm}}

Let $(\nu,\lambda,\nu^\infty)$ be an admissible measure-valued solution with initial data $v_0$. By Proposition \ref{claim1} it suffices to generate the lifted admissible measure-valued subsolution by a suitable sequence of subsolutions. By abuse of notation, let us still denote the lifted Young measure by $(\nu,\lambda,\nu^\infty)$, and its energy by $E(t)$. 

By Theorem \ref{YMapprox} we obtain a sequence of smooth subsolutions $\bar{w}^k=(\bar{v}^k,\bar{u}^k)$ uniformly bounded in $L^{\infty}_t(L^2_x\times L^1_x)$ and a sequence of discrete oscillation measures $\nu^k$ that are piecewise constant with respect to $\mathcal{Q}_k$, such that 
\eqref{approxconvergence}, \eqref{approxenergy} and \eqref{e:approxinitial} hold.

Let us fix $Q\in\mathcal{Q}_k$ for the moment. By Proposition \ref{discretegeneration} there exists a sequence 
$$
\tilde w^{k,n}=(\tilde v^{k,n},\tilde u^{k,n})\in C_c^\infty(Q;\Real^d\times\mathcal{S}^d_0)
$$
of subsolutions
such that 
$$
\|\tilde w^{k,n}\|_{L^{\infty}(Q)}\leq \langle\nu^k|_Q,e\rangle
$$ 
and generating $\nu^k$ on $Q$. 
Next, for $l\geq 1$, let $\bar{w}^{k,l}=(\bar{v}^{k,l},\bar{u}^{k,l})$ be defined as the average of $(\bar{v}^k,\bar{u}^k)$ over cubes of size $2^{-l}\frac{1}{k}$, 
i.e.
$$
(\bar{v}^{k,l},\bar{u}^{k,l})(x,t):=\fint_{Q_i^l}(\bar{v}^k,\bar{u}^k)\quad\textrm{ for $(x,t)\in Q_i^l$ for all $Q_i^l\in\mathcal{Q}_{2^lk}$.}
$$
Since $(\bar{v}^k,\bar{u}^k)\in L^{\infty}_t(L^2_x\times L^1_x)\cap C^{\infty}(\Real^d\times[0,T])$, we have 
$$
\sup_t\int_{\Real^d}|\bar{v}^{k,l}-\bar{v}^k|^2\,dx+\sup_t\int_{\Real^d}|\bar{u}^{k,l}-\bar{u}^k|\,dx\to 0\quad\textrm{ as }l\to\infty.
$$
Consequently, from \eqref{e:epointwise} we deduce 
$$
\sup_t\biggl|\int_{Q_t}e(\bar{w}^k(x,t)+\tilde w^{k,n}(x,t))\,dx-\int_{Q_t}e(\bar{w}^{k,l}(x,t)+\tilde w^{k,n}(x,t))\,dx\biggr|\to 0
$$
as $l\to\infty$, where $Q_t=\{x:\,(x,t)\in Q\}$, and the convergence is uniform in $k,n$, owing to the $L^{\infty}$ bound on $\tilde w^{k,n}$ as well as the energy bound \eqref{approxenergy} on
$\nu^k$.
On the other hand, since $e\in C(\Real^d\times\mathcal{S}^d_0)$ is convex, from \eqref{convexestimate} we have  
$$
\limsup_{n\to\infty}\sup_t\int_{Q_t}e(\bar{w}^{k,l}+\tilde w^{k,n}(x,t))\,dx\leq \langle\mathcal{T}_{\bar{w}^{k,l}}\nu^k|_{Q},e\rangle\,.
$$
Thus, given any $\delta>0$ we can first choose $l$ and then $n$ sufficiently large, so that
$$
\sup_t\int_{Q_t}e(\bar{w}^{k}+\tilde w^{k,n}(x,t))\,dx\leq \langle\mathcal{T}_{\bar{w}^{k}}\nu^k|_{Q},e\rangle+\delta\,.
$$

Summing over all cubes in $\mathcal{Q}_k$ and taking a suitable diagonal subsequence, we obtain $(\tilde v^{k},\tilde u^k)$ such that $(\bar{v}^k+\tilde{v}^k,\bar{u}^k+\tilde{u}^k)$ is uniformly bounded in $L^{\infty}_t(L^2_x\times L^1_x)$ (this follows from \eqref{approxenergy} and \eqref{e:uniformbd}), consists of subsolutions, and
$$
(\bar{v}^k+\tilde{v}^k,\bar{u}^k+\tilde{u}^k)\overset{\mathbf{Y_{2,1}}}{\longrightarrow} (\nu,\lambda,\nu^\infty)\,.
$$
Furthermore, 
$$
\sup_t\int_{\Real^d}e(\bar{v}^k+\tilde{v}^k,\bar{u}^k+\tilde{u}^k)\,dx\leq \textrm{esssup}_tE(t)+\frac{1}{k}
$$
and, by \eqref{e:approxinitial}, 
$$
\int_{\Real^d}|\bar{v}^k(x,0)+\tilde{v}^k(x,0)-v_0(x)|^2\,dx=\int_{\Real^d}|\bar{v}^k(x,0)-v_0(x)|^2\,dx\to 0
$$
as $k\to\infty$. Therefore the sequence $(\bar{v}^k+\tilde{v}^k,\bar{u}^k+\tilde{u}^k)$ satisfies the conditions in Proposition \ref{claim1}b). Proposition \ref{almostenergythm} follows,
since $(\nu,\lambda,\nu^{\infty})$ is admissible and therefore $\textrm{esssup}_tE(t)=\tfrac{1}{2}\int_{\Real^d}|v_0|^2\,dx$. 
 \qed
 
\subsubsection{Proof of Theorem \ref{energythm}}\label{proofenergythm}

As above, let $(\nu,\lambda,\nu^{\infty})$ be an admissible measure-valued
solution, and, by abuse of notation, let us denote by
$(\nu,\lambda,\nu^{\infty})$ also the corresponding lifted Young measure on
$\Real^d\times\mathcal{S}_0^d$ as in Section
\ref{lifting}. In the proof of Proposition \ref{almostenergythm} above, we showed that
for every $k\in\mathbb{N}$ there exists a subsolution
$(v^{k},u^k)\in C^{\infty}(\Real^d\times[0,T])$
such that
\begin{align}
(v^k,u^k)\overset{\mathbf{Y_{2,1}}}{\longrightarrow} (\nu,\lambda,\nu^{\infty})\textrm{ as }k\to\infty\\
\left|\frac{1}{2}\int|v^{k}(x,0)|^2dx-\frac{1}{2}\int|v_0(x)|^2dx\right|<\tfrac{1}{k}\label{epsinitial}\\
\sup_{t\in[0,T]}\int_{\Real^d}
e(v^{k},u^{k})dx<\frac{1}{2}\int|v_0|^2dx+\tfrac{1}{k}\label{epsenergy}
\end{align}
In this final step,
we want to deduce from this the full statement of Theorem
\ref{energythm}. We will do so by using the argument from
\cite{euler2} for the construction of wild initial data. A slight modification
of the proof of Proposition 5.1 in \cite{euler2} yields the following statement:
\begin{prop}\label{p:wild}
Let $(v,u)\in L^{\infty}_t(L^2_x\times L^1_x)\cap C^{\infty}(\Real^d\times[0,T])$
be a subsolution, and let 
$$
\bar{e}\in C(\Real^d\times[0,T])\cap
C([0,T];L^1(\Real^d))
$$ 
such that 
$e(v(x,t),u(x,t))<\bar{e}(x,t)$ for all
$(x,t)\in\Real^d\times[0,T]$. Then there exists a sequence of 
subsolutions $(v^n,u^n)\in
C^{\infty}(\Real^d\times(0,T])$ such that
\begin{equation*}
v^n\to v\quad\textrm{ in }C([0,T];L^2_w(\Real^d)),
\end{equation*}
\begin{equation*}
e(v^n(x,t),u^n(x,t))<\bar{e}(x,t)\quad\forall\,(x,t)\in\Real^d\times(0,T]
\end{equation*}
and
\begin{equation}\label{e:wildinitial}
\frac{1}{2}|v^n(x,0)|^2=\bar{e}(x,0)\quad\textrm{ for a.e.~$x\in\Real^d$.} 
\end{equation}
\end{prop}
Indeed, the proof follows along the lines of the proof
of Proposition 5.1 in \cite{euler2}, by iterating the Claim from there 
with the functional
$\int_{\Omega_0}\tfrac{1}{2}|v(x,0)|^2-\bar{e}(x,0)\,dx$ instead of
$\int_{\Omega_0}\tfrac{1}{2}|v(x,0)|^2-1\,dx$, and only considering the 
functions for times $t\geq 0$.

With this assertion at hand,
we proceed as in the proof of Proposition \ref{claim1}. Choose
$\bar{e}_k\in C(\Real^d\times[0,T])\cap
C([0,T];L^1(\Real^d))$ so that
\begin{itemize}
\item $e(v^k(x,t),u^k(x,t))<\bar{e}_k(x,t)$ for all $(x,t)\in \Real^d\times[0,T]$;
\item $\sup_t\int_{\Real^d}\bar{e}_k-e(v^k,u^k)\,dx<\tfrac{1}{k}$;
\item $\int_{\Real^d}\bar{e}_k(x,t)\,dx\leq\int_{\Real^d}\bar{e}_k(x,0)\,dx$ for all $t\in [0,T]$.
\end{itemize}
Such a choice is possible because
of \eqref{epsinitial}-\eqref{epsenergy}.

Then, apply Proposition \ref{p:wild} to each $(v^k,u^k)$ with $\bar{e}_k$. We obtain sequences 
$(v^{k,n},u^{k,n})$, and as in the proof of Proposition \ref{claim1} we can extract a diagonal subsequence $n=n(k)$
such that
\begin{align*}
v^{k,n(k)}-v^k\to 0\textrm{ in }L^2_{loc}(\Real^d\times(0,T))\,.
\end{align*}
Moreover, we have that $v^{k,n(k)}(t=0)\rightarrow v_0$ strongly in $L^2$ if $n(k)$ is chosen sufficiently large: Indeed, on the one hand, we know that $v^{k,n(k)}(t=0)\rightharpoonup v^k(t=0)$ weakly (cf. Proposition \ref{p:wild}). On the other hand, by (\ref{e:wildinitial}), the choice of $\bar{e}_k$, (\ref{epsenergy}), and (\ref{epsinitial}), we have
\begin{equation*}
\norm{v^{n,k}(t=0)}^2_{2}-\norm{v^k(t=0)}^2_2<\frac{8}{k}.
\end{equation*}
These two facts imply that $v^{k,n}(t=0)\rightarrow v^k(t=0)$ strongly as $n\rightarrow\infty$, and hence, if $n(k)$ is suitably chosen, (\ref{epsinitial}) yields $v^{k,n(k)}(t=0)\rightarrow v_0$ strongly in $L^2$ as $k\rightarrow\infty$.

Finally, with each $v^{k,n(k)}$ we argue as in the proof of Proposition \ref{claim1}b) to
find exact Euler solutions close to $v^{k,n(k)}$ in $C([0,T];L^2_w(\Real^d))$ and hence, due to the choice of
$\bar{e}^k$, in $L^2_{loc}(\Real^d\times(0,T))$. Therefore these Euler solutions generate the same
Young measure. Furthermore, their initial values $v^{k,n(k)}(t=0)$ are close to $v_0$, and because $v^{k,n(k)}$ satisfies \eqref{e:wildinitial}, they are admissible. 
This proves Theorem \ref{energythm}. \qed

\subsubsection{Proof of Corollary \ref{c:density}}

As observed by DiPerna and Majda in \cite{dipernamajda}, given any initial data $v_0\in L^2(\Real^d)$, a sequence of Leray solutions with viscosity $\to 0$ generates a measure-valued solution of Euler. Moreover, such a measure-valued solution has initial data $v_0$, and it is easy to see from the energy bound that  it will be admissible. Hence, as in the previous Subsection \ref{proofenergythm}, we find a sequence of subsolutions $(v^{k,n(k)},u^{k,n(k)})$ such that (passing to a subsequence if necessary)
\begin{equation*}
\norm{v^{k,n(k)}(t=0)-v_0}_2<\frac{1}{k},
\end{equation*}
\begin{equation*}
\frac{1}{2}|v^{k,n(k)}(x,0)|^2=\bar{e}_k(x,0)\hspace{0.3cm}\text{for a.e. $x\in\Real^d$},
\end{equation*}
and
\begin{equation*}
e\left(v^{k,n(k)}(x,t),u^{k,n(k)}(x,t)\right)<\bar{e}_k(x,t)\hspace{0.3cm}\text{for a.e. $x\in\Real^d$ and all $t>0$}.
\end{equation*}

Hence, by virtue of Theorem \ref{dense} (which is Proposition 3.3. in \cite{euler2}), for each such $(v^{k,n(k)},u^{k,n(k)})$ there exist infinitely many solutions for Euler with initial data $v^{k,n(k)}(t=0)$ and energy density $\bar{e}_k$ for all times $t\geq0$. By choice of $\bar{e}_k$ (see Subsection \ref{proofenergythm}), these solutions are admissible. This shows that in every $L^2$-neighbourhood of $v_0$ there exists initial data admitting infinitely many admissible solutions of Euler and thus the corollary is proved.
\qed

Finally we remark that if we choose $\bar{e}_k$ such that $\int\bar{e}_kdx$ is constant in time, then the Euler solutions
obtained in this way will conserve energy. This shows the existence of
{\em energy-conserving} weak solutions for a dense subset of initial data.

%\bibliography{General}

\end{document}